\pdfoutput=1 
\documentclass[11pt, final]{elsarticle}
\usepackage[left=1.5cm,top=1.5cm,right=1.5cm,bottom=1.5cm]{geometry}
\usepackage{booktabs}
\usepackage{amsmath}
\usepackage{algorithm}
\usepackage{algpseudocode}
\usepackage{xcolor}
\usepackage{tikz}
\usepackage{multirow}
\usepackage{lscape}
\usepackage{color}
\usepackage{rotating}
\usepackage{pdflscape}
\usepackage{todonotes}
\usepackage{caption}
\usepackage{longtable}
\usepackage{xr-hyper}
\usepackage[draft,commentmarkup=footnote]{changes}
\usepackage{natbib}
\usepackage[hidelinks]{hyperref}
\usepackage{adjustbox}
\usetikzlibrary{spy}
\allowdisplaybreaks 

\usepackage{cleveref} 
\usepackage{subfig}
\usepackage{pgfplots}
\usetikzlibrary{patterns}
\usepackage{csvsimple}
\usepackage{pgfplotstable}

\usepgfplotslibrary{colorbrewer}
\pgfplotsset{compat = 1.15, cycle list/Set1-8}
\usetikzlibrary{pgfplots.statistics, pgfplots.colorbrewer}

\makeatletter
\newcommand*{\addFileDependency}[1]{
  \typeout{(#1)}
  \@addtofilelist{#1}
  \IfFileExists{#1}{}{\typeout{No file #1.}}
}
\makeatother

\newcommand*{\myexternaldocument}[1]{%
    \externaldocument{#1}%
    \addFileDependency{#1.tex}%
    \addFileDependency{#1.aux}%
}

\myexternaldocument{Supplement_Appendix}

\definechangesauthor[name={Jake}, color=red]{jw}
\definechangesauthor[name={Andreas}, color=magenta]{AE}
\definechangesauthor[name={Yuan}, color=blue]{YS}

\usepackage{algorithmicx,eqparbox,array}

\algrenewcommand{\algorithmiccomment}[1]{\hfill// \eqparbox{COMMENT\thealgorithm}{#1}}
\algnewcommand{\LongComment}[1]{\hfill// \begin{minipage}[t]{\eqboxwidth{COMMENT\thealgorithm}}#1\strut\end{minipage}}

\usepackage{titlesec}
\titleformat{\paragraph}
{\normalfont\normalsize\itshape}{\theparagraph}{1em}{}
\titlespacing*{\paragraph}
{0pt}{3.25ex plus 1ex minus .2ex}{1.5ex plus .2ex}

\biboptions{authoryear}

\usepackage{etoolbox}
\patchcmd{\pprintMaketitle}
 {\ifvoid\absbox\else\unvbox\absbox\par\vskip10pt\fi}
 {\ifvoid\absbox\else\clearpage\unvbox\absbox\par\vskip30pt\fi}
 {}{}
\patchcmd{\pprintMaketitle}
 {\hrule\vskip12pt}
 {}
 {}{}
\patchcmd{\pprintMaketitle}
 {\hrule\vskip12pt}
 {}
 {}{}
\appto{\pprintMaketitle}{\clearpage}

\begin{document}

\begin{frontmatter}

\title{Ranking Constraint Relaxations for Mixed Integer Programs Using a Machine Learning Approach}


\author{Jake Weiner\fnref{myfootnote}$^{a}$}
\author{Andreas T.~Ernst$^{b}$} 
\author{Xiaodong Li$^{a}$}
\author{Yuan Sun$^{b}$}

\address{$^{a}$ School of Computing Technologies, RMIT University, 124 La Trobe Street, Melbourne, Australia}
\address{$^{b}$  School of Mathematics, Clayton Campus, Monash University, Melbourne, Australia}

\fntext[myfootnote]{Corresponding Author.\\ 
Email Addresses: s3730771@student.rmit.edu.au (J. Weiner), andreas.ernst@monash.edu (A. Ernst), xiaodong.li@rmit.edu.au (X. Li), yuan.sun@monash.edu (Y. Sun)}

\begin{abstract}
Solving large-scale Mixed Integer Programs (MIP) can be difficult without advanced algorithms such as decomposition based techniques. Even if a decomposition technique might be appropriate, there are still many possible decompositions for any large MIP and it may not be obvious which will be the most effective. This paper presents a comprehensive analysis of the predictive capabilities of a Machine Learning ranking (ML) function for predicting the quality of Mixed Integer Programming (MIP) decompositions created via constraint relaxation. In this analysis, the role of instance similarity and ML prediction quality is explored, as well as the benchmarking of a ML ranking function against existing heuristic functions. For this analysis, a new dataset consisting of over 40000 unique decompositions sampled from across 24 instances from the MIPLIB2017 library has been established. These decompostions have been created by both a greedy relaxation algorithm as well as a population based multi-objective algorithm, which has previously been shown to produce high quality decompositions. In this paper, we demonstrate that a ML ranking function is able to provide state-of-the-art predictions when benchmarked against existing heuristic ranking functions. Additionally, we demonstrate that by only considering a small set of features related to the relaxed constraints in each decomposition, a ML ranking function is still able to be competitive with heuristic techniques. Such a finding is promising for future constraint relaxation approaches, as these features can be used to guide decomposition creation. Finally, we highlight where a ML ranking function would be beneficial in a decomposition creation framework.
\end{abstract}

\begin{keyword}
Machine Learning \sep  Integer Programming \sep Heuristics \sep Large Scale Optimization
\end{keyword}

\end{frontmatter}

\section{Introduction}

As Mixed Integer Programming (MIP) problems continue to increase in both size and complexity, new solution techniques are required to keep up with the increased problem complexities. This continual increase is reflected in the popular MIP benchmark dataset MIPLIB, where from 2003 to 2017 the largest instances in the dataset increased from just over 200,000 variables to over 1,400,000 variables \citep{Achterberg2006, Gleixner2021}. As most traditional methods require non-polynomial time complexities for solving these MIPs, this increase in problem size cannot be addressed by improvements made in computer processing power alone. One popular approach used to solve these large scale problems is via decomposition, a long standing and powerful idea for which three prominent methods -  Benders Decomposition (BD) \citep{Geoffrion1972}, Dantzig-Wolfe Reformulation (DWR) \citep{Vanderbeck2006} and Lagrangian Relaxation (LR) \citep{Geoffrion1974} are most commonly used. All methods are often able to generate tighter bounds than linear programming based relaxation bounds, as well as feasible solutions if a repair heuristic is used or if the method is embedded within a Branch and Bound algorithm. As a result, these decomposition techniques are still prominent in the research community today, despite being decades old. 

Whilst problem decomposition has been extremely successful, this approach has only been applicable to a limited number of problem types in the literature, for which the problem structure lends itself to be nicely decomposable. Ideally, when practitioners are attempting to solve a new and previously unseen MIP without a known decomposable structure, it would be beneficial if an automatic decomposition framework could be used to determine if the problem might benefit from a decomposition based approach and which decomposition is most promising. An automatic decomposition framework such as this might expose numerous problems to be decomposable when otherwise a decomposition approach would not be attempted. At the same time, it would allow practitioners to try a decomposition approach without having significant domain specific knowledge about the problem. 

Whilst an automatic decomposition framework might be extremely beneficial, until recently, there has been only limited work carried out in this area. Three studies related to the creation of an automatic decomposition framework were carried out in \citep{Bergner2015a}, \citep{Khaniyev2018} and in \citep{Weiner2020}, as well as an automated DWR tool called GCG \citep{Gamrath2010}. In \citep{Bergner2015a}, the authors represent the original constraint matrix as a hypergraph, after which a $k$-way partitioning algorithm is run in an attempt to create $k$ equal sized independent blocks. The decompositions are then solved within a DWR framework. The only user input required for this procedure is the number of pre-determined blocks into which the problem should be decomposed and a number of dummy nodes for partition balancing purposes. One significant drawback for this approach is the a-priori requirement of knowing the correct number of blocks into which the problem should decompose, a feature which can vary significantly on an instance by instance basis. In \citep{Khaniyev2018} and \citep{Weiner2020}, the authors address this issue and propose automatic decomposition methods without requiring any user input. In \citep{Khaniyev2018}, this is achieved via a community detection algorithm and a uniquely defined decomposition quality metric labelled as the Goodness score. This score is derived from both non-zero percentages in the subproblems and the proportion of relaxed constraints in the border. The authors in \citep{Weiner2020} instead treat the problem of decomposition as a multi-objective one, using the Non-Dominated Sorting Genetic Algorithm II (NSGA-II) algorithm \citep{Deb2002} to automatically create decompositions which minimize the number of constraints relaxed and the size of the largest subproblem simultaneously. What all three papers potentially lack is that they only consider a small set of features to guide the automatic decomposition process, with the number of blocks, the number of relaxed constraints and the size of the subproblems primarily used as features. In addition, for population based heuristics such as the NSGA-II algorithm described in \citep{Weiner2020}, when multiple decompositions are available, having a single objective ranking function becomes important.


More recently, a new branch of work has emerged, considering numerous features and using a Machine Learning (ML) method to identify both when a decomposition approach should be attempted and what decomposition is most promising amongst several candidates. The authors in \citep{Kruber2017} train a variety of classifiers to answer the question of whether or not a decomposition is likely to result in producing better bounds than simply using a generic solver on the original problem. Whilst some success was shown using this approach, this line of research only helps to address whether or not a decomposition approach may be appropriate. It does not, however, provide insight as to what decomposition should be selected for solving when considering multiple candidates. The authors in \citep{Basso2020} are the first to attempt to fill in this gap regarding the ranking of decomposition quality a-priori based on a variety of decomposition features. Whilst some promising results were found, ultimately this line of research is still in its infancy, with only limited success being demonstrated. In particular, we believe further analysis is required regarding 1) What role instance similarity plays in ranking performance; 2) if a ML ranking function is able to accurately predict decomposition quality in test instances not sampled in the training set; 3) how a ML ranking function performs when benchmarked against other heuristic measures of decomposition quality; and 4) are there any constraint features which can be identified as being important to decomposition quality, as this could be used in future sampling procedures.

In this paper, we provide an extensive analysis regarding the effectiveness of a ML approach to rank decomposition quality, where a combination of solve time and bound quality are used as a metric for quantifying decomposition quality. To do so, we have created a significant dataset of over 40,000 unique decompositions sampled from 24 instances in the MIPLIB2017 library. These decompositions were created using the multi-objective approach as described in \citep{Weiner2020} which is able to produce a rich set of Pareto optimal decompositions using both border and subproblem metrics, as well as a greedy random selection algorithm to ensure the dataset maintains some diversity. We show the effect of instance similarity on ML performance, and demonstrate how a ML ranking function can produce state-of-the-art predictions when benchmarked against four established heuristic functions found in both the literature and in online solvers. We also show which constraint features could be beneficial in guiding future automatic decomposition frameworks involving constraint relaxation. Finally, we provide a summary and possible future work on how our findings may be used in a heuristic to create fully automated decompositions.

The rest of this paper is structured as follows. Section \ref{section:Background} introduces the background and related work. Our main approach is described in Section \ref{section:Approach}. Section~\ref{section:Experimental Setup} details the experimental design, Section~\ref{section:Results} presents the results and Section~\ref{Section:Conclusion} then concludes the paper.

\section{Background}\label{section:Background}
In this section we briefly describe the Lagrangian Relaxation framework for empirical evaluation of a decomposition's effectiveness and previous works which attempt to quantify decomposition quality through either heuristic or ML based methods, the greedy and NSGA-II frameworks used to create decompositions as described in \citep{Weiner2020}.

\subsection{Lagrangian Relaxation}\label{section:LR}

Lagrangian Relaxation (LR) is a popular decomposition technique that is commonly implemented in the Operations Research (OR) community. LR is most often applied when the `right' set of complicating constraints are identified and relaxed, thereby decomposing the original problem into multiple independent subproblems which are significantly easier to solve. Problems suitable for this type of decomposition display angular matrix patterns such as those seen in Figure~\ref{fig:angular_matrix}. In this angular constraint matrix, often referred to as a Bordered Block-Diagonal matrix structure, there are complicating constraints $A_c$ and independent block structures $D_1,\ldots,D_k$. If these complicating constraints were removed, the subproblems $D_1,\ldots,D_k$ are naturally able to be decomposed and solved independently, with the hope that solving these subproblems to optimality can be achieved significantly faster than the original problem as a whole, although there is no guarantee.


\begin{figure}[tb!]
\centering
$\begin{bmatrix}
\textit{D}_1 &  &  &  &  \\
 & \textit{D}_2 &  &  & \\
 &  &  \textit{D}_3  & & \\
&  &  & \ddots & \\
&  &  &   & \textit{D}_k\\
\cdots & \cdots & \textit{A}_c & \cdots & \cdots
\end{bmatrix}$
\quad
\caption{Bordered Block-Diagonal Matrix Structure. Decomposed subproblems $D_1,\ldots,D_k$ contain independent variables and can be solved independently if the complicating constraints $A_c$ are removed.}
\label{fig:angular_matrix} 
\end{figure}

The new Lagrangian objective function is formed by shifting the relaxed constraints to the objective function, with a penalty term (Lagrangian Multiplier) attached. A basic implementation of the primal and Lagrangian dual formulations are show in Eqns~\eqref{eqn:LR simple1}-\eqref{eqn:LR simple2} respectively, where $x$ are the decision variables, $c$ are the associated costs, $A$ and $b$ are the constraint matrix and resource constraints and $\lambda$ are the Lagrangian Multipliers introduced. Whilst this formulation is representative of a Binary Linear Program, it can easily be extended to a Mixed Integer Program. Solving this new Lagrangian dual with the optimal Lagrangian Multipliers is often able to provide high quality problem bounds. In addition, solutions to the new Lagrangian dual problem are often almost feasible for the original primal problem, and can often be made feasible through an appropriate repair heuristic, or when embedded in a Branch and Bound framework.
\begin{align} 
    \max\limits_{x} \quad f(x) = &\, c^Tx \quad  \text{s.t.}\  \   Ax \leq b,\quad D x \leq d,\quad x \in \{0,1\}^n  \label{eqn:LR simple1} \\ 
    \min\limits_{\lambda\ge0} \;LR(\lambda) =&\max \left\{ c^Tx + \lambda^T\, (b - Ax) : Dx\le d,\ x\in\{0,1\}^n\right\}\label{eqn:LR simple2}
\end{align}
The choice of decomposition here corresponds to the choice of how the constraints are partitioned between the $A$ and $D$ matrix. The effectiveness of the LR approach is influenced by the difficulty in solving the non-smooth dual problem (minimisation over $\lambda$) and the trade-off between the best lower bound that can be achieved (tight relaxation) and the ease of evaluating $LR(\lambda)$ (small, easily solvable subproblems). As both bound quality and solve time are important metrics in evaluating the effectiveness of a LR based approach, for this paper we assign scores to decompositions using both of these metrics as described in Section~\ref{Sect: Training and Testing}.

To find the optimal Lagrangian Multipliers, there exist a wide variety of potential approaches, including Subgradient Optimization \citep{Fisher2004}, Bundle Methods \citep{Barahona2000a}, Coordinate Ascent \citep{Wedelin1995a} and Hybrid techniques \citep{Ernst2012a,Weiner2021}. For the purposes of this paper we will avoid the question of the solution approach by evaluating $LR(\lambda)$ only once with the optimal dual values $\lambda^*$ from the corresponding constraints in the LP relaxation of the primal problem (\ref{eqn:LR simple1}). Warm starting the LR procedure with the optimal dual LP values also provides a clear metric as to the usefulness of the decomposition. If the Integrality Property holds, the optimal LP solution is equivalent to the Lagrangian dual $LR(\lambda^*)$ and is indicative that the decomposition is not useful when solved using LR. This is unless optimising the Lagrangian dual is able to be achieved faster than solving the LP, which is very unlikely to occur when solving general MIPs without a specialised algorithm. This warm start approach for finding the initial, or potentially final Lagrangian Multipliers is also proposed in \citep{Geoffrion1974}. As even only solving the Lagrangian objective function can be computationally expensive for large problems, the aim is to find some quality metrics that accurately predict performance based on directly observable attributes of the decomposition.

\subsection{Heuristic Decomposition Quality Metrics}

Within the literature there exists only a limited number of studies investigating the quality of decompositions for general MIPs. These approaches can be split into heuristic based methods and Machine Learning based approaches. The heuristic based approaches are presented in \citep{Bergner2015a, Khaniyev2018, Weiner2020} and investigate decomposition qualities using defined heuristic measures. Two of these heuristics  \citep{Bergner2015a, Khaniyev2018} are used as benchmarks in Section~\ref{section:Results} and the final approach \citep{Weiner2020} is used to generate the set of decompositions for training and testing purposes.

\paragraph{Relative Border Area}
In 2015, the first paper \citep{Bergner2015a} demonstrating the potential for an automatic decomposition approach was presented, in the context of solving decompositions via a Dantzig-Wolfe Reformulation process. To do so, the authors represent the constraint matrices of general Mixed Integer Programs as hypergraphs, for which they then solve the minimum weight balanced $k$-partition problem \citep{Karypis1999} in an attempt to create $k$ equal sized subproblems whilst reducing the number of relaxed constraints required to find such a partition. Without presenting a framework to discover which single decomposition should be solved via DWR, the authors instead perform an exploratory analysis on different parameter selections, in particular the number of $k$ subproblems created. The authors also present a heuristic measure which they suggest can be used to find good decompositions from amongst several candidates, referring to this measure as the Relative Border Area (RBA). Formally, the RBA is presented as: $\frac{m_l \times n + m\times n_l - m_l \times n_l}{m \times n}$ where $m_l$ is the number of linking constraints, $n_l$ is the number of linking variables, $m$ is the total number of constraints and $n$ is the total number of variables. Small RBA values are indicative of high quality decompositions. The authors note however, that for 26 out of 39 instances tested, decompositions without linking variables performed the best. When no linking variables are considered, the RBA is simply $\frac{m_l}{m}$, the percentage of constraints relaxed. Whilst this is a fairly simple heuristic, the authors note `DWR unfolds its potential when relaxing a few linking constraints is enough to decompose the problem into more tractable subproblems'.


\paragraph{Goodness Score}
Another automatic decomposition framework is presented in \citep{Khaniyev2018} and uses a community detection algorithm to maximise a uniquely defined decomposition quality metric the authors have defined as the Goodness score. As noted by the authors, a good decomposition can be calculated by analysing both the subproblem and border components separately, which are labelled as $Q$ and $P$ measures respectively. These metrics are then combined to give an overall Goodness score. The subproblem and border scores can be calculated as follows:
\begin{description}
\item[Subproblem Components]  With regards to subproblem components, three factors are identified in being important to decomposition quality:
\begin{enumerate}
    \item Granularity - A decomposition with a large number of subproblems is desirable, as in theory, solving smaller subproblems is computationally beneficial.
    \item Homogeneity - The subproblems should ideally be equally sized. Homogeneity is especially important when solving subproblems on parallel processors, ensuring that the largest subproblem does not dominate the solve time taken. 
    \item Isomorphism - Ideally, subproblems should be identical not only in size, but in objective function coefficients and right-hand side values. This results in a single subproblem required to be solved at every iteration of the Lagrangian Relaxation algorithm instead of all subproblems. 
\end{enumerate}
The overall proxy used to measure these subproblem statistics is presented as $Q = (\sum_{i=1}^{K}\frac{nz_i}{m}[1 - \frac{nz_i}{m}])$ where $nz_i$ is the number of nonzero entries in block $i: \; \forall i \in K$ blocks and $m$ is the number of non-zeroes in the block-diagonal component of the constraint matrix.
\item[Border Component]
Following on from the results presented in \citep{Bergner2015a}, the quality of the border component calculated within the Goodness score is the percentage of constraints relaxed, albeit in a simple exponential decay function, as the authors note that there is non-linearity in the correlation between border size and optimality gap. The overall proxy used to measure the quality of the border is presented as $P = (e^{-\lambda(\frac{b}{M})})$ where b is the number of constraints in the border and M is the total number of constraints.
\end{description}
The final Goodness score for a decomposition is calculated as $Q \times P$ and is bounded between $[0,1]$. Decompositions which have larger Goodness scores are considered to be higher quality.

\paragraph{Multi-objective Approach}\label{section:NSGA-II}
A multi-objective approach for automatically generating good decompositions is presented in \citep{Weiner2020} and uses the well known Non Dominated Sorting Genetic Algorithm II (NSGA II) to create and evolve a population of decompositions. In a manner similar to the Goodness score, two objectives are minimised which have been shown to result in high quality decompositions. The multi-objective algorithm aims to minimise both the number of constraints relaxed (small border area) and the size of the largest subproblem for each decomposition. This approach was used to generate both the training and test datasets used in this paper, as this framework is able to generate a large number of good decompositions in a relatively short amount of time.

\paragraph{Heuristic Measures outside of the literature}

Within the GCG solver \citep{Gamrath2010} a quality metric for general decomposition with linking constraints only can be found in the open source code. For benchmarking, we have labelled this heuristic as GCG Open Source (GCG OS), and is calculated as: $(0.6 \times \frac{m}{M}) + 0.01 + (0.2 \times (1 - \min\{d_1,\ldots,d_k\}))$ where $m$ is the number of constraints relaxed, $M$ is the total number of constraints and $d_k$ is the non-zero density of the coefficient matrix for subproblem $k$. 

Another heuristic that exists is the Max-White (MW) score presented on the MIPLIB 2017 \citep{Gleixner2021} website and is used to determine if an instance is suitable for decomposition or not. The MW score is calculated as $1 - (\frac{s + t}{nvars * ncons})$ where $t = nvars_1 \times ncons_1$ and $s = \sum_{i=2}^{nb+1} nvars_i * ncons_i$ where $nb$ is the number of subproblems, with the first subset being the border.

\subsection{Random Sampling and Machine Learning}

To our knowledge, the authors in \citep{Basso2020} are the first to address the question of how to rank decomposition quality a-priori based on instance and decomposition features. To do so, a pseudo random sampling algorithm was used to generate decompositions, and both classification and regression models were trained to predict Pareto optimal solutions (where bound quality and solution time form the two objectives), and then subsequently rank the selected decompositions according to distances from the closest Pareto optimal solution. This framework is referred to as the data-driven process. For classification, in initial experiments when decompositions from each base instance were used as the test set and the classifier was trained on decompositions from all other instances, poor results were observed from all but 5 out of the 36 instances tested. Furthermore, even when decompositions from all instances were included in the training set, the test precision score ($P$), where $P = \frac{TP}{TP + FP}$, is only 0.0714, with a significant number of false positives predicted. The authors give two potential reasons for the poor results that arose during classification. Firstly, the authors note that positive decompositions in different instances might not have enough common characteristics, and suggest that perhaps more positive decompositions are required in the dataset. Because of this finding, we have included decompositions in our dataset which have been created via the NSGA-II algorithm as described in \citep{Weiner2020}, as these decompositions have been shown to be of greater quality than those produced by a pseudo-random sampling process. A second reason given for poor classification results is that instances in the dataset are likely to be structurally different, and therefore including decompositions from all instances for training could be detrimental if good decomposition patterns vary depending on the structure of the instance. The authors suggest further research is needed to investigate the effect of instance similarity and the performance of a trained classifier. As such, in this paper we focus on problems that include a network structure in the hope of increasing instance similarity. We provide an analysis on the importance of instance similarity in Section~\ref{results:intersubproblem_performance}. The next phase of the data-driven approach was to train a regressor in order to then rank the positive decompositions chosen by the classification step. Whilst there are some interesting preliminary results, as noted by the same authors in \citep{Basso2018}, the performance of the regressor tended to be quite poor when used as a ranking function. In \citep{Basso2018}, the authors attempt to define a new ranking function based on dominance percentage instead of Pareto distances, with some improvements being noted. Finally, the benchmarking carried out in \citep{Basso2020} only compares a limited number of decompositions created via a data-driven process against decompositions created by the static detectors in GCG. What is lacking, is how a ML ranking function compares to other heuristic based decomposition ranking functions, such as the Max-White score \citep{Gleixner2021}, the decomposition score for bordered block-diagonal decompositions as found in GCG's Open Source code (GCG OS), the Relative Border Area (RBA) metric discussed in \citep{Bergner2015a} and the Goodness Score presented in \citep{Khaniyev2018}.




\section{Approach}\label{section:Approach}

The approach carried out in this paper consists of four main tasks that are described in more detail below:
\begin{enumerate}
    \item Decomposition generation using both greedy-random and multi-objective approaches.
    \item Decomposition  post-processing, including the removal of both redundant constraints and duplicate decompositions.
    \item Establishing MIP bounds by solving the Lagrangian function for decompositions using the optimal dual values from the LP relaxation as Lagrangian Multipliers. 
    \item Analysis of decomposition results. The analysis carried out addressed three main areas:
    \begin{itemize}
        \item Investigate the relationship between instance similarity and prediction performance.
        \item Benchmarking of ML methods against Heuristic techniques for both similar and dissimilar instances.
        \item Investigate the relationship between Relaxed Constraint features and prediction quality.
    \end{itemize}
\end{enumerate}


\subsection{Decomposition Generation} 

As noted by the authors in \citep{Basso2020}, using a randomised sampling algorithm from an arbitrary MIP is unlikely to yield promising decompositions. The authors therefore implemented a pseudo random sampling algorithm, selecting constraints with a probability in proportion to their sparsity, in a manner similar to the greedy-random decompositions as discussed in \citep{Weiner2020}. For the dataset used in this paper, we used decompositions created both by the NSGA-II algorithm and greedy-random approaches as described in \citep{Weiner2020}. We included the NSGA-II generated decompositions as these were previously shown to produce high quality decompositions when compared to greedy-random generated decompositions. Additionally, in order to make sure the dataset was not too biased towards decompositions that only consider the two metrics used in the multi-objective approach, the greedy-random approach was also included. For both approaches, the constraint matrix of the MIP being sampled is translated into a hypergraph, where rows are represented by hyperedges and columns form the nodes within the hypergraph. Once the hypergraph is created, the set of constraints to relax are selected and consequently the corresponding hyperedges are removed from the hypergraph. A Breadth First Search (BFS) algorithm is then run to identify the independent partitions which now exist, representing the independent subproblems created. The BFS algorithm is able to run in $\mathcal O(nz)$ where $nz$ is the number of non-zero entries in the constraint coefficient matrix. In this paper we have also introduced two post processing steps to remove unnecessary decomposition components in order to better filter the dataset.

\paragraph{Greedy-random Decompositions}
A greedy-random approach was used to generate a diverse range of decompositions, with a bias in relaxing constraints with a larger number of non-zeroes, in the hope that doing so would lead to more subproblems which are smaller. To create a decomposition, constraints are sorted according to the number of variables they contain, then they are iterated over with a probability of being relaxed equal to $p_i = \frac{|V_i|}{S_V}\ \times Q \times |C|$, where $p_i$ is the probability of constraint $i$ being relaxed, $V_i$ is the set of variables contained in constraint $i$, $S_V$ the sum of all variable counts across the set of all constraints $C$ in the original problem and $Q$ is the desired proportion of constraints to be relaxed. This iterative loop is run until the desired proportion of relaxed constraints has been selected.

\paragraph{NSGA-II Decompositions}
The NSGA-II decompositions were created as described in \citep{Weiner2020}.  In the NSGA-II implementation, the fitness function consists of minimising two objectives: 1) the number of constraints relaxed and 2) the size of the largest subproblem (measured by the number of variables). The initial population is seeded with some greedy-random solutions, containing a variety of different percentages of relaxed constraints, in an effort to assist with search exploration. An arbitrary number of individuals are initialised with varying numbers of constraints relaxed, from 5\% to 99\% of the total number of constraints.

\subsection{Decomposition Post-processing}


Two decomposition post-processing steps have been created to eliminate both redundant relaxed constraints and duplicate decompositions. For the purposes of this paper, when referring to a constraint as redundant, it is done so with regards to how likely it is to affect to either the solve time or bound quality of the decomposition if this constraint is relaxed. It is not suggesting that removing such a constraint does not change the set of feasible solutions. Post decomposition creation, we have identified two ways in which a relaxed constraint may be considered redundant and not provide any meaningful contribution to decomposition quality.

Examining the relaxed constraints and the subproblems which occur as a result, a relaxed constraint is considered redundant if either:

\begin{enumerate}
    \item All variables associated with the constraint are already a subset of the variables in one of the subproblems.
    \item All variables associated with the constraint are \textbf{only} found in single variable subproblems.
\end{enumerate}

\begin{figure}[!t]
\centering
\includegraphics[width=9cm]{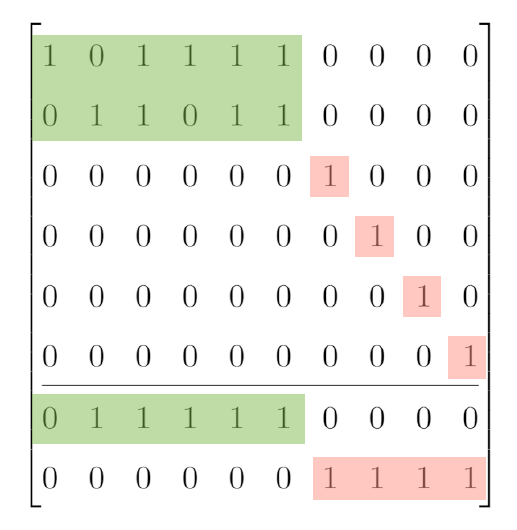}\vspace{-2mm}
\caption{Constraint Redundancy Processing. Shown in green and red are two examples of relaxed constraints which can be considered to be redundant. Shown in green is an example of a redundant constraint in which all non-zeroes belong to a subproblem, in which case there is no decomposition of the subproblem when this constraint is relaxed. Shown in Red is a constraint in which all non-zeroes are found only in single variable subproblems, for which when solved provide no additional tightening of the LP bounds.}
\label{fig:Redundant Constraints}
\end{figure}

In the first definition of a redundant constraint, if all variables of a relaxed constraint are contained within a subproblem, moving the constraint from the border back to the subproblem is unlikely to affect the solve time of the subproblem in any substantial way, however it can result in tighter bounds.  As such we have deemed this constraint to be redundant for the purposes of decomposition. In the second definition of a redundant constraint, if all variables in a constraint are found in only single variable subproblems, the solution of the variables in these subproblems is simply set to the variable bounds, providing no tightening of the LP bound. An example of two constraints which are considered redundant is shown in Figure~\ref{fig:Redundant Constraints}. Shown in green is a constraint which is considered redundant, as all non-zeroes contained within are already part of an existing subproblem. Shown in red is a constraint in which all non-zeroes are contained \textbf{only} in single variable subproblems. These single variable subproblems can arise as a result of constraint relaxation, where variables are no longer part of any constraints and are simply constrained by their bounds. Finally, after all constraint post-processing, all duplicate decompositions are removed from the dataset. 


\subsection{Establishing MIP Bounds}

To establish MIP bounds, the LP relaxation of the original problem was solved first, with the optimal dual values $\lambda^*$ used as the Lagrangian Multipliers. The Lagrangian Multipliers were then used to solve the Lagrangian function. Each independent subproblem was solved to within 1\% of optimality, as at this stage only bound quality is considered and there can be significant time spent by modern solvers trying to prove optimality. This evaluation of the Lagrangian function was run in an effort to close the integrality gap \citep{Fisher2004}, improving upon the LP relaxation bound and providing a proxy of the decomposition bound quality.

Solving the Lagrangian function with Lagrangian Multipliers set to the optimal dual values from the corresponding constraints in the LP relaxation of the primal problem was chosen for several reasons and is also recommended in \citep{Geoffrion1974}. First, due to the significant number of decompositions created ($\geq 40,000$), the computing time required to run the full Lagrangian Relaxation algorithm with a multiplier update procedure for several iterations would be computationally prohibitive. As each subproblem is attempted to be solved to optimality, this process can be extremely computationally expensive. Secondly, the introduction of a multiplier update procedure adds another layer of complexity and could result in decomposition qualities being reliant on the update procedure, for which there exist a variety of potential methods available, such as Subgradient Optimization \citep{Fisher2004}, Bundle Methods \citep{Barahona2000a}, Coordinate Ascent \citep{Wedelin1995a} and Hybrid techniques \citep{Ernst2012a,Weiner2021}. In addition, for some of these methods, a stochastic element exists for which multiple runs would need to be carried out. 


Due to the relaxation of constraints, there may exist subproblems which contain only single variables. For these subproblems, they can be solved to optimality using the variable bounds instead of via solving a MIP. After relaxation, subproblems were sorted in ascending order by variable sizes and given a runtime limit proportional to the square number of variables, as a linear relationship between subproblem variable sizes and MIP solution times seems unlikely. The solve time limit for each subproblem was calculated as $t_i = (v^2_{i} / \sum_{i = 1}^{K} {v^2_{i}}) \times CPU $ where $t_i$ is the subproblem run time, $v_i$ is the number of variables in subproblem $i$, $K$ is the set of subproblems and $CPU$ is the total CPU budget allocated to solve the Lagrangian Objective function with the optimal LP dual values. For the final and largest subproblem, the full remaining run-time was allocated.

If a subproblem was unable to be solved to optimality within the given run-time limit, the best bound found during the Branch and Bound process was used as the subproblem solution. If the root node was unable to be processed in this time, an additional $60$(s) of CPU runtime was given in an attempt to establish a valid bound for the subproblem and to avoid potential bound errors found by the MIP solver used. In the rare event that the root node was still unable to be processed, the subproblem was solved as a LP and the LP bound was used, with no additional time taken to solve the subproblem as a LP counted towards the total runtime.


\subsection{Dataset Selection}

In addressing the first part of the analysis process as described in Section~\ref{section:Approach} regarding how instance similarity affects prediction performance, we selected three different network problem types from the MIPLIB2017 library \citep{Gleixner2021}, forming what we refer to as the Network dataset. These three problem types - Network Design (ND), Fixed Cost Network Flow (FCNF) and Supply Network Planning (SNP) each contain instances generated from the same optimisation model but using different problem data. In order to address how well a ML approach would then generalise to randomly selected instances, we also selected 10 instances from the MIPLIB2017 library containing properties which represent a broad range of potential unseen instances. We refer to the dataset comprised of these instances as the Random (Rand) dataset. The 10 instances from the MIPLIB library were chosen such that 1) A significant proportion of discrete variables in order to facilitate tightening of the LP bound, although, this is not necessarily always a requirement, as even a few binary or integer variables can still significantly tighten the LP bound; 2) A reasonable number of non-zeroes in order for constraint processing and decomposition creation were able to run within a reasonable amount of time; 3) A pre-processing step to ensure that the LP bounds for most of the instances was not optimal, a necessary requirement in order to effectively rank decompositions; 4) Some instances which are easily solved to optimality quickly, perhaps indicating that relaxing no constraints is considered the best decomposition. The full set of instances and instance statistics is shown in Table~\ref{tab:Instance Statistics}.


\begin{table*}[!t]
\setlength\tabcolsep{10.5pt}
\def\arraystretch{1}
\caption{Instance Statistics: Shown for instance in the dataset are the total number of variables (var), number of binary (bin), integer (int) and continuous (cont) variables, number of constraints (constr) and the number of non-zeroes (nz). Each instance is also shown which problem type it belongs to.}
\label{tab:Instance Statistics}
\begin{tabular}{lrrrrrrr}
\toprule
Instance Name & Var & Bin & Int & Cont & Constr & NZ & Problem Type \\ \midrule
cost266-UUE & 4161 & 171 & 0 & 3990 & 1446 & 12312 & ND \\
dfn-bwin-DBE & 3285 & 2475 & 0 & 810 & 235 & 9855 & ND \\
germany50-UUM & 6971 & 0 & 88 & 6883 & 2088 & 20737 & ND \\
ta1-UUM & 2288 & 0 & 605 & 1683 & 439 & 5654 & ND \\
ta2-UUE & 9241 & 1188 & 0 & 8053 & 2687 & 26533 & ND \\ \midrule
g200x740 & 1480 & 740 & 0 & 740 & 940 & 2960 & FCNF \\
h50x2450 & 4900 & 2450 & 0 & 2450 & 2549 & 12152 & FCNF \\
h80x6320d & 12640 & 6320 & 0 & 6320 & 6558 & 31521 & FCNF \\
k16x240b & 480 & 240 & 0 & 240 & 256 & 960 & FCNF \\ \midrule
snp-02-004-104 & 228350 & 167 & 167 & 228016 & 126512 & 463941 & SNP \\
snp-04-052-052 & 221438 & 4546 & 4546 & 212346 & 129662 & 459205 & SNP \\
snp-06-004-052 & 328461 & 494 & 494 & 327473 & 183168 & 668716 & SNP \\
snp-10-004-052 & 538777 & 815 & 815 & 537147 & 300348 & 1097780 & SNP \\
snp-10-052-052 & 549021 & 11059 & 11059 & 526903 & 320836 & 1138760 & SNP \\ \midrule
splice1k1 & 3253 & 3252 & 1 & 0 & 6505 & 1761020 & Rand \\
neos-4954672-berkel & 1533 & 630 & 0 & 903 & 1848 & 8007 & Rand \\
dws008-01 & 11096 & 6608 & 0 & 4488 & 6064 & 56400 & Rand \\
traininstance2 & 12890 & 5278 & 2602 & 5010 & 15603 & 41531 & Rand \\
neos-4338804-snowy & 1344 & 1260 & 42 & 42 & 1701 & 6342 & Rand \\
neos-4387871-tavua & 4004 & 2000 & 0 & 2004 & 4554 & 23496 & Rand \\
30n20b8 & 18380 & 18318 & 62 & 0 & 576 & 109706 & Rand \\
air05 & 7195 & 7195 & 0 & 0 & 426 & 52121 & Rand \\
blp-ic98 & 13640 & 13550 & 0 & 90 & 717 & 191947 & Rand \\
air03 & 10757 & 10757 & 0 & 0 & 124 & 91028 & Rand \\ \bottomrule
\end{tabular}
\end{table*}





\subsubsection{Features Considered}
Whilst there is no consensus amongst the optimisation community as to exactly which instance and decomposition features directly correspond to decomposition quality, it has been found that minimising the number of linking constraints and maximising the number of similarly sized subproblems can close the integrality gap and result in faster decomposition runtimes respectively \citep{Bergner2015a, Khaniyev2018}. Minimising the number of relaxed constraints and the size of the largest subproblem were also shown to be beneficial in generating good decompositions \citep{Weiner2020}. In order to exploit the Machine Learning model approach, we use a more extensive list of features, relying on both subproblem statistics and border statistics. The full list of features can be found in Table~A1 in the Appendix and contain many of the features used in \citep{Kruber2017} and \citep{Basso2020}, for which the authors note there exists combination of features which are important to decomposition quality.

Shown in Appendix Figures~A1 and A2 are the spread of the normalised feature data for both the Network and Random datasets respectively. Within these datasets, the percentage of decompositions in which all features are between $[Q1 - 1.5 \times IQR, Q3 + 1.5\times IQR]$, where $Q1$ is the 25th percentile, $Q3$ is the 75th percentile and $IQR$ is the interquartile range, are 9.68\% and 4.43\% respectively. Such a spread highlights the potential difficulty of using a ML approach, as feature data can vary significantly between instances, containing a significant number of `outlier' points.




\subsubsection{Training and Testing}\label{Sect: Training and Testing}

Finally, for the purposes of this paper, the best decomposition for each instance is the one that achieves the lowest score, which we have defined as a combination of both bound quality (represented as the optimality gap ( $|\frac{UB-LB}{UB}| * 100$) and solve time. Using the Weighted Sum Model \citep{Triantaphyllou2000}, we attribute equal importance to both gap values $(g)$ and solve times ($t$), $score = 0.5*g + 0.5*t$, although in future works this weighting can change depending on user requirements. All decompositions are assigned scores on an instance by instance basis, where all optimality gaps and solve times are normalised using Min-Max normalisation. A Min-Max normalisation process was also used to normalise all features not already in a 0-1 scale range. A variety of regression models were then trained to predict decomposition score, in an effort to rank the quality of decompositions a-priori based on the features considered. Because a ML model is only useful if decomposition quality can be predicted for unseen instances, an approach approximating Leave-One-Out-Cross-Validation \citep{Sammut2010} method was used for testing and is further described in Section~\ref{section:Results}. The regression models trained include some of the most popular linear and non linear regression models as well as two different ensemble methods including:
\begin{enumerate}
    \item Linear Regression with Ridge Regularisation (Ridge)
    \item Linear Regression with Lasso Regularisation (Lasso)
    \item Support Vector Regression (SVR) with a radial basis function kernel
    \item K Nearest Neighbours Regression (KNN)
    \item Random Forest Regressor (RF)
    \item Multi-Layer Perceptron (MLP)
    \item Stacking Ensemble using Ridge, Lasso, SVR, KNN and MLP trained regressors and a Linear Regression Estimator (Stacking)
    \item Voting Ensemble using Ridge, Lasso, SVR, KNN and MLP trained regressors (Voting)
\end{enumerate}

\section{Experimental Setup}\label{section:Experimental Setup}
Hypergraph partitioning and NSGA-II algorithms were run on an Intel i7-7500U CPU and all LR tests were carried out on the Multi-modal Australian ScienceS Imaging and Visualisation Environment (MASSIVE) network, which runs on an Intel Xeon CPU E5-2680 v3 processor. CPLEX 12.8.0 using a single thread was used to solve all MIP subproblems and LP benchmarks. All LR runs were given a limited runtime of 300 CPU seconds for feasibility reasons, as the total number of decompositions tested was in excess of 40,000 as shown in Table~\ref{tab:Decomposition Counts}. The NSGA-II algorithm was run using the Pagmo framework \citep{Biscani2019} with default parameter settings. These settings include: Crossover Probability = 0.95, Distribution index for Crossover = 10.0, Mutation Probability = 0.01, Distribution index for Mutation = 50.0. The NSGA-II algorithm was run using 300 generations with population sizes of 32. Whilst good convergence was demonstrated for smaller generation numbers, using the percentage of Pareto optimal solutions in generation $n$ which were also found in generation $n-1$ as a performance measure, due to the fast run-time of our hypergraph partitioning algorithm, we included 300 generations in order to attain a richer and more diverse dataset. It should be noted that for some of the Supply Network Planning instances, the full NSGA-II algorithm was unable to be completed within the allocated runtime. This is due to the implementation of the NSGA-II algorithm, for which it was discovered post analysis there was a slight logic error within the crossover implementation. For future work, this error within the crossover implementation may be addressed to speed up the search time, however for our purposes a sufficient number of decompositions were generated for analysis. For the greedy-random decompositions, 999 decompositions were generated for each instance, with 111 decompositions created for specified relaxed constraint percentages ranging from 10\% to 90\% of total constraints. Included amongst the decompositions tested was also a decomposition in which no constraints were relaxed. While over 10,000 decompositions for each instance were generated with the above approaches, after removal of duplicates the number of decompositions available for training and testing is in the range of 500--2,500 per instance as shown in Table~\ref{tab:Decomposition Counts}.

All ML models were trained with default parameter settings as presented in the Scikit-learn library \citep{scikit-learn} with Python 3.6.13. The only exception to this were the alpha regularisation parameters chosen for the Lasso and Ridge Linear Regression models, for which a range of regularisation parameters was tested, as the training time for these models was relatively insignificant. For the Lasso and Ridge models, the alpha parameters selected were 0.001 and 0.01 respectively. 

\section{Results}\label{section:Results}

Multiple experiments were carried out to identify 1) The effect of instance similarity on the prediction quality of the ML methods; 2) How does a ML ranking approach compare to other heuristic based ranking functions found in the literature and in open source solvers?; 3) Which relaxed constraints features are important to decomposition quality? and 4) How well do ML methods perform when tested on randomly selected instances from the MIPLIB library with seemingly no similarities in problem structure? When presenting the scores of the best decomposition selected by each of the ranking methods, this is the best score from amongst the top 8 decompositions identified by the ranking method. A similar approach is used in \citep{Basso2020}, as in practice it is fairly trivial to solve 8 decompositions in parallel on most modern computer architectures and therefore the best decomposition amongst the top 8 selected can be easily be identified. Finally, when presenting the scores from the best decompositions selected by the different ranking methods, an additional Min-Max normalisation was carried out for decomposition scores on an instance by instance basis. This additional normalisation more easily shows the quality of the selected decompositions by the ranking methods, as a score of 0 indicates the best decomposition from the population was selected and a score of 1 indicates the worst decomposition from the population was selected.

\subsection{Instance Similarity and Performance} \label{results:intersubproblem_performance}

As noted in \citep{Basso2020}, instance similarity is a potential reason for the relatively poor results found via the classification experiments the authors carried out when predicting Pareto optimal decompositions. Therefore, we investigate the significance of instance similarity on prediction quality by comparing ML models trained and tested on decompositions from the same problem type and models trained and tested on decompositions from other problem types. For each instance and model type, we trained the model on all decompositions from the same problem type \textbf{excluding} the test instance, an approach approximating Leave-One-Out-Cross-Validation (LOOCV). The predictions of this model was then compared against the same model trained on all decompositions from each of the other problem types.

\begin{figure*}[ht]
\centering
\subfloat[~Testing Problem Type: Network Design]{\includegraphics[scale=.5]{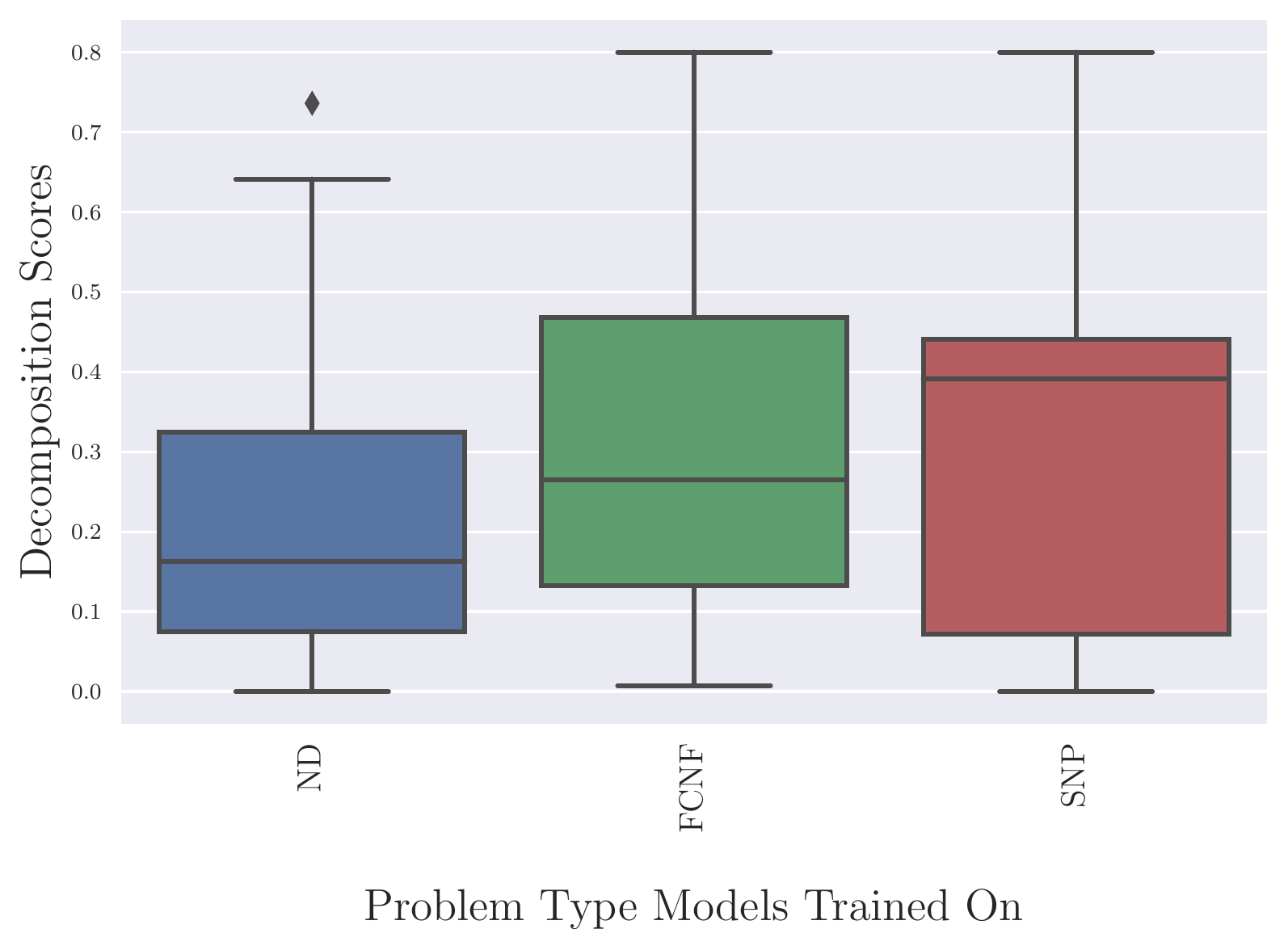}}\,
\subfloat[~Testing Problem Type: Fixed Cost Network Flow]{\includegraphics[scale=.5]{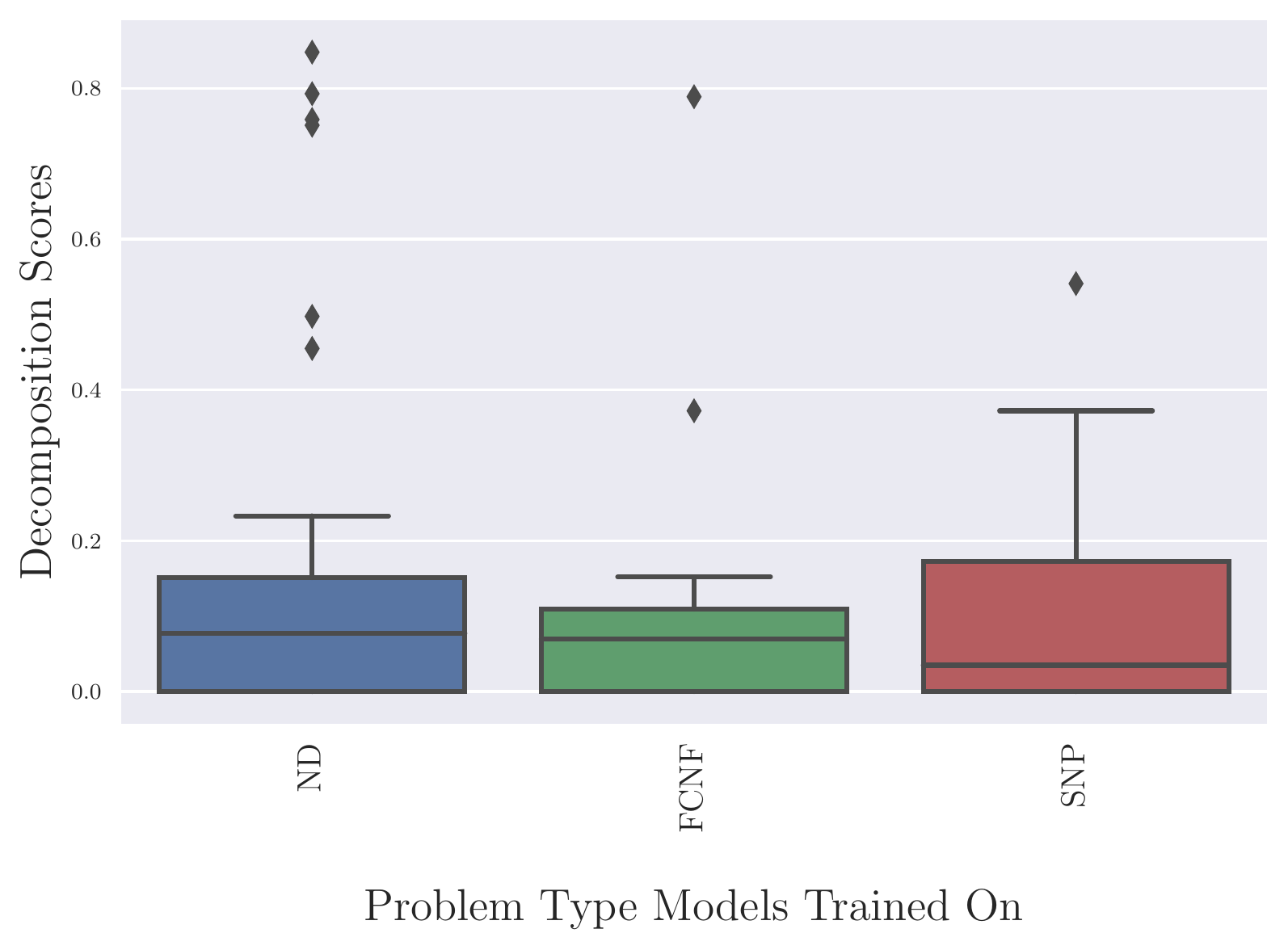}}
\\
\subfloat[~Testing Problem Type: Supply Network Planning]{\includegraphics[scale=.5]{images/Intersubproblem/fixed_cost_network_flow.pdf}}
\vspace{0.1 cm}  
 \caption{Instance Similarity and Prediction Quality: Shown for when instances from each problem type in the Network dataset are used as test instances (unseen by the trained models), are the decomposition scores of the best predicted decompositions selected by each of the ML models when trained on the different problem types. For models trained on the same problem type as the test problem type, the model is trained on all instances in the problem type except for the test instance. When the training problem type is different to the testing problem type, all instances in the training problem type are used for training. The decomposition scores range from 0 (the best decomposition was selected) to 1 (the worst decomposition was selected).}
\label{fig: Intersubproblem_boxplot}
\end{figure*}

The decomposition scores of the best decompositions selected by the ML models are presented in Tables~A4, A5 and A6 in the Appendix, with boxplots of the results shown in Figure~\ref{fig: Intersubproblem_boxplot}. These results show that models trained and tested on instances from the same problem type are able to better predict decomposition quality than models trained and tested on different problem. To detect if there is any statistical significance between models trained on instances from the same problem type as opposed to instances trained on instances from a different problem type, we carried out pairwise comparisons using the Friedman Aligned Rankings and a post-hoc analysis as described in \citep{Derrac2011}. This non-parametric statistical test was chosen as there is no underlying assumption that the predicted decomposition scores would follow a normal distribution. The Friedman Aligned Rankings and subsequent post-hoc analysis were used as the number of comparison methods is relatively low ($=3$) and as such this statistical approach is recommended \citep{Derrac2011}. The z-scores for each pairwise comparison can be calculated as $z = (\tilde{R_{i}} - \tilde{R_{j}})/ \sqrt{\dfrac{k(n + 1)}{6}}$, where $\tilde{R_{i}}$ and $\tilde{R_{j}}$ are the average Friedman Aligned Rank for the control algorithm and comparison algorithm respectively, $k$ is the number of comparison algorithms and $n$ is the number of samples. From the calculated z-scores, an unadjusted p-value can be found from the table of normal distribution $N(0,1)$. In every pairwise comparison (except for one) between models trained and tested on instances from the same problem type, and models trained and tested on instances from a different problem type, there is a statistical significance detected ($p \leq 0.05$) as shown in Tables~A4, A5 and A6. Based on these results we conclude that instance similarity does play an important part in the prediction capabilities of a ML based ranking function. Even amongst instances which all contain a  Network structure, a ML ranking function clearly performs better when test instances are from the same problem type as the training instances, albeit from different data sources.

\begin{figure}[!t]
\centering
\includegraphics[width=0.9\linewidth]{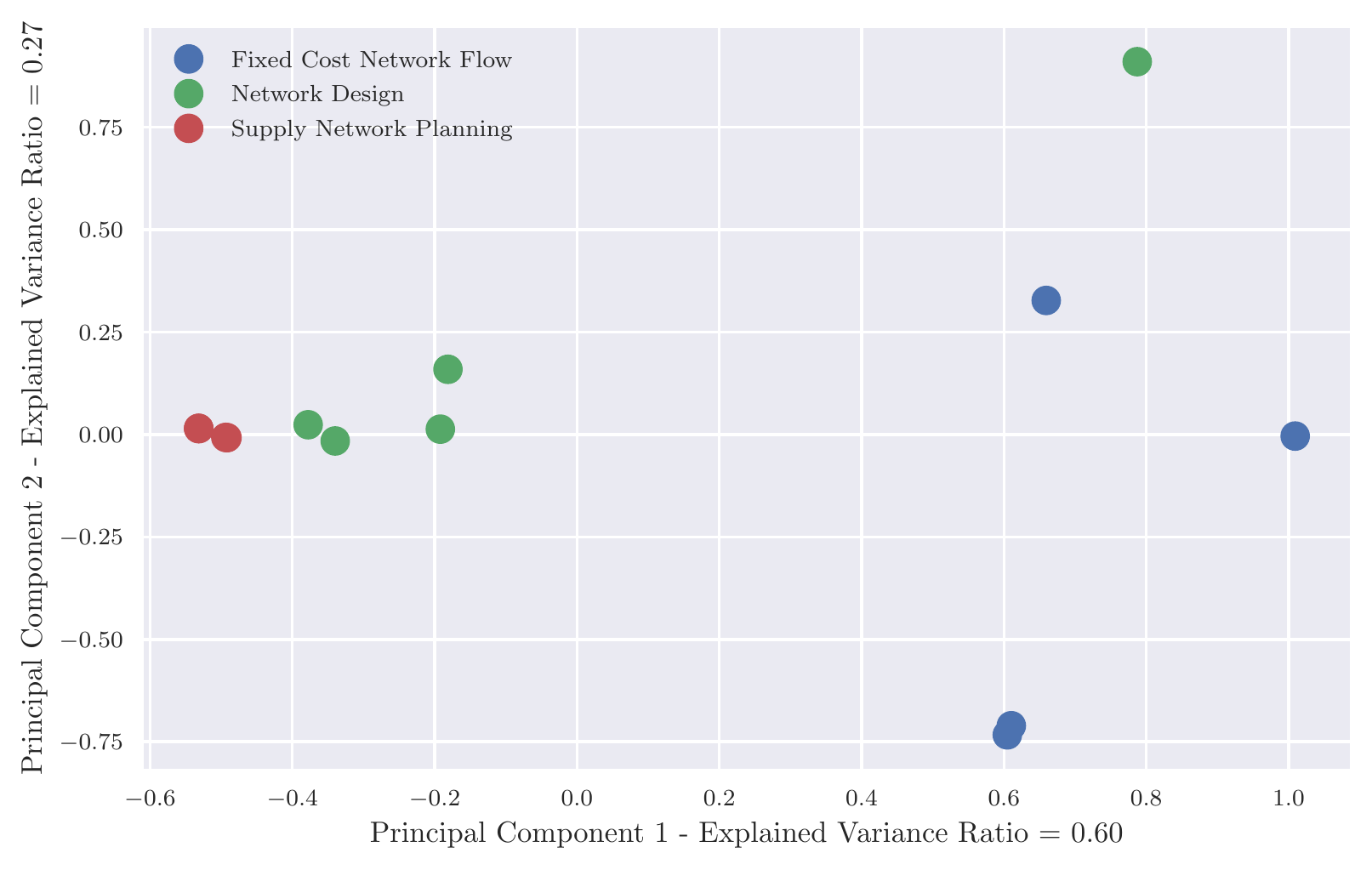}\vspace{-2mm}
\caption{Principal Component Analysis of the Network Design, Fixed Cost Network Flow and Supply Networking instances using the features described in Table~A7. Shown in this figure are the two principle components, comprising 87\% of the explained variance. As can be seen, instances from the Supply Network Planning and Network Design problem types display good clustering tendencies, indicating high similarity between the instances within these problem types. The Fixed Cost Network Flow problem type shows poor clustering, indicating instances within this problem type are less similar to each other than the other Network problem types.}
\label{fig:PCA Network Problems}
\end{figure}

In order to visualise if instance similarity can be captured only by considering some relatively simple features, we carried out a Principle Component Analysis (PCA) using a linear dimensionality reduction on all instances from the three network problem types, using the instance derived features as shown in Table~A7. A kernel based PCA was also carried out, including third degree polynomial, radial basis function, sigmoid and cosine kernels, however these showed no improvements upon instance separation and clustering tendencies. As can be seen in Figure~\ref{fig:PCA Network Problems}, the Supply Network Planning problem type displays excellent clustering tendencies, indicating high similarity between the instances within this problem type. Similarly, the Network Design instances also seem to be more similar to one another in this feature space, except for one outlier instance. For the Fixed Cost Network Flow problem type, the instances seem to be less similar to each other. These PCA visualisations appear to be in line with our findings in Figure~\ref{fig: Intersubproblem_boxplot}, for which the models both trained and tested on the Network Design and Supply Network Planning problem types seem to perform significantly better than when the models are trained and tested on different problem types.

An additional observation can be made regarding instance similarity in this feature space we have explored. In this feature space, the Network Design and Supply Network Planning instances seem to be more similar to each other than instances from the Fixed Cost Network Flow problem type. From this, we could assume models trained and tested on these instances might perform better than models trained on either of these problem types and tested on the Fixed Cost Network Flow problem type, or vice versa. From our findings this is not always the case, indicating that there are still other features which could be important in determining instance similarity.

\subsection{State-of-the-art Benchmarking}

\begin{figure*}[!htb]
\centering
\includegraphics[width=0.7\linewidth]{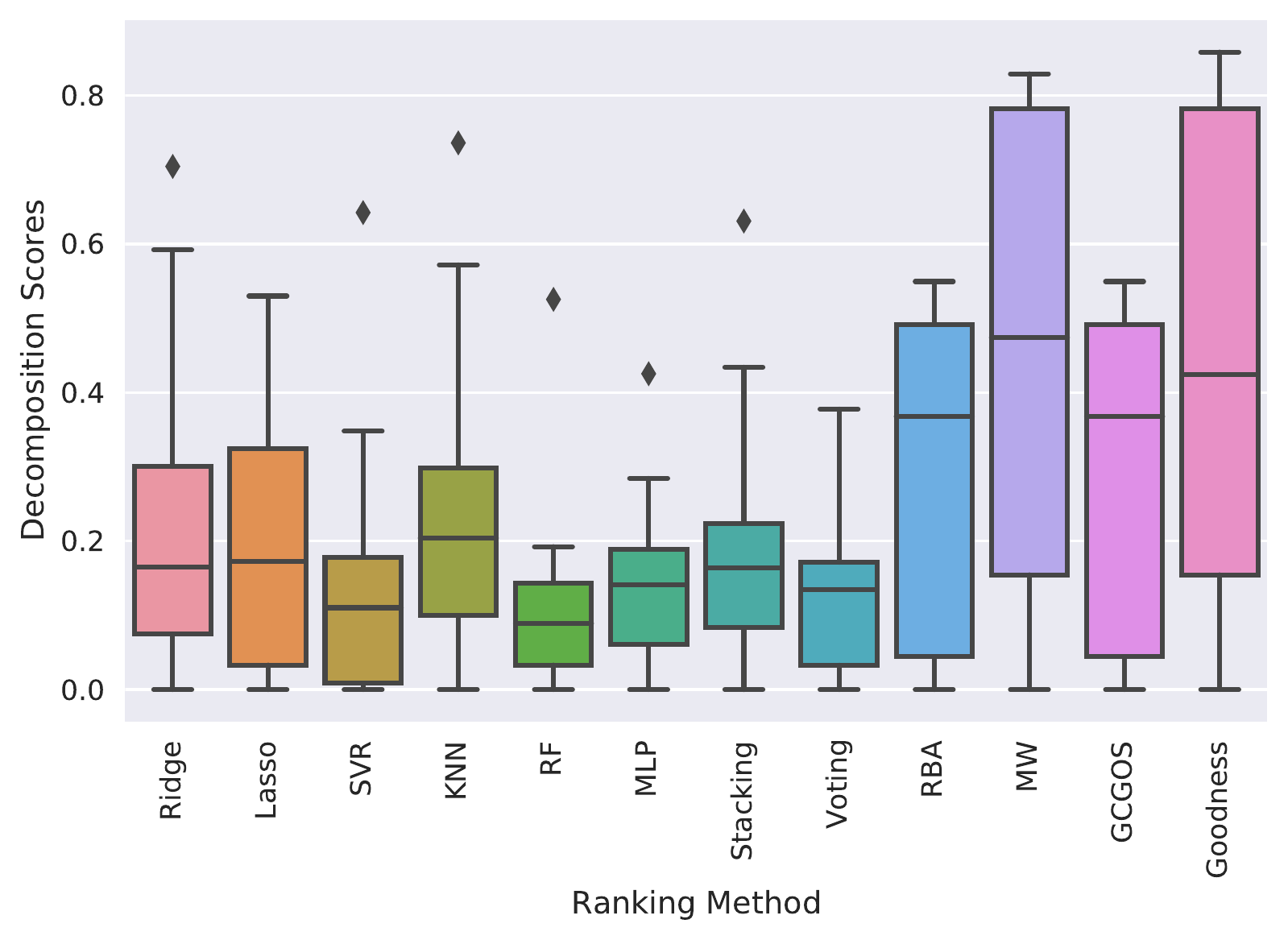}\vspace{-2mm}
\caption{Prediction Benchmarking. Shown for each ranking method is a boxplot of the scores for the best decompositions selected by each ranking method across the 14 test instances in the Network dataset. The decomposition scores range from 0 (the best decomposition in the test set) to 1 (the worst decomposition in the test set).}
\label{fig:network_prediction_boxplot}
\end{figure*}

\begin{figure*}[!htb]
\centering
\includegraphics[width=0.7\linewidth]{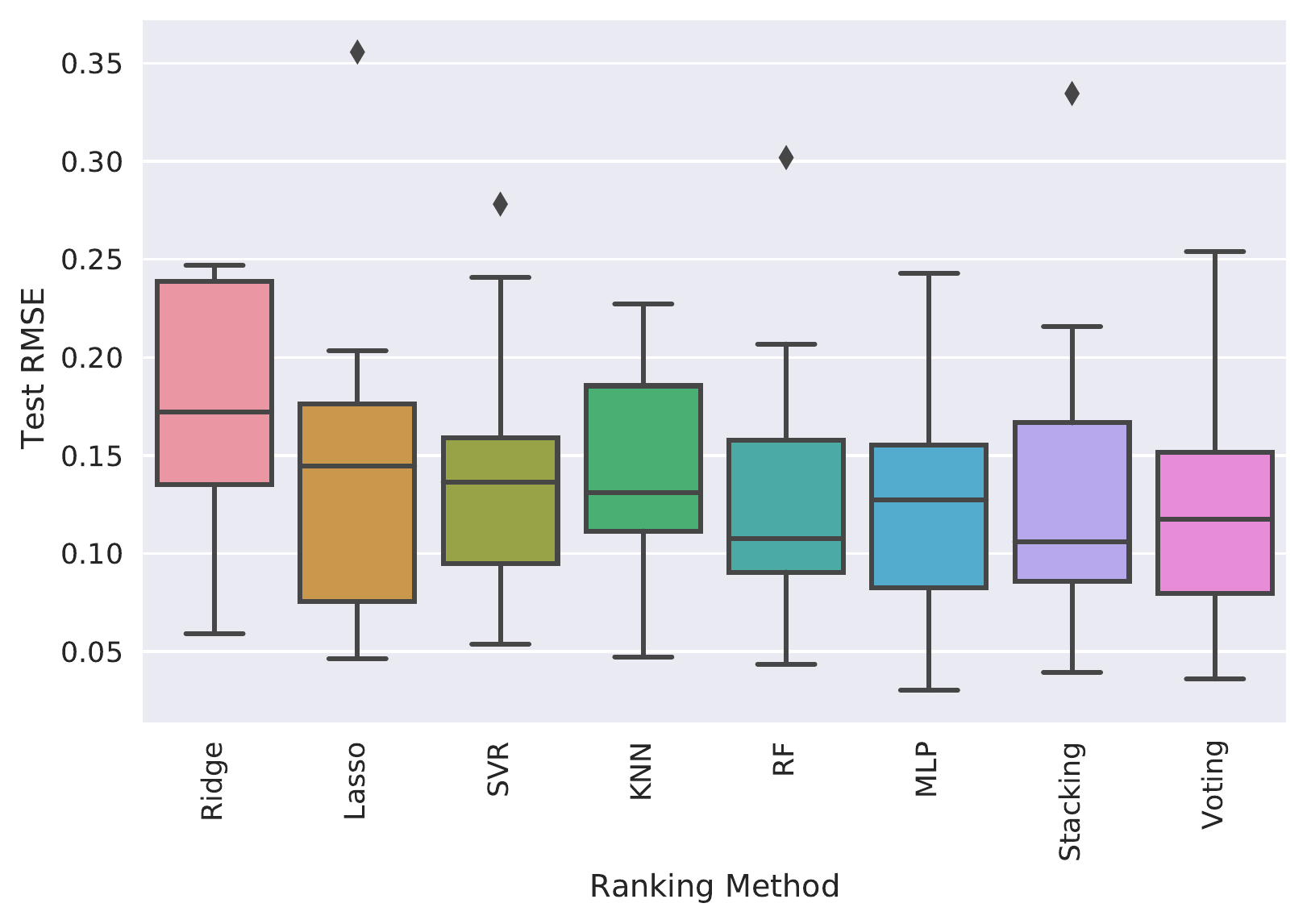}\vspace{-2mm}
\caption{ML RMSE Scores. Shown for each ML ranking method is a boxplot of the test RMSE scores for the 14 test instances in the Network dataset.}
\label{fig:network_rmse_boxplot}
\end{figure*}

In order to investigate how effective a ML model is at selecting the best quality decomposition from a population of decompositions, we benchmarked the various trained ML models against four heuristic approaches which have been used for measuring decomposition quality, all of which can be found in the literature and in solver open source code. These include the decomposition score metric that can be found in GCG's Open Source codebase (GCG OS) \citep{Gamrath2010}, the Goodness score as presented in \citep{Khaniyev2018} (Goodness), the Max White score as described on MIPLIB 2017 (MW) \citep{Gleixner2021} and the Relative Border Area metric as described in \citep{Bergner2015a}. Shown in Table~\ref{tab:RankingMethodsRaw} are the decomposition scores of the best decompositions selected by each ranking method. For the ML methods, each model was trained on all decompositions from the Network dataset \textbf{excluding} the test instance, an approach approximating Leave-One-Out-Cross-Validation. These results are visualised in Figure~\ref{fig:network_prediction_boxplot}, with the test Root Mean Square Errors (RMSE) for the ML predictions shown in Figure~\ref{fig:network_rmse_boxplot}.


\begin{landscape}
\begin{table*}[!p]
\centering
\caption{Decomposition Prediction Results: Shown for each ranking method is the decomposition score for the best decomposition selected by the method. The decomposition scores range from 0 (the best decomposition in the test set) to 1 (the worst decomposition in the test set). Shown in bold are the scores for the best predicted decomposition for each instance. Also shown is the average decomposition score predicted by each ranking method.}
\label{tab:RankingMethodsRaw}
\begin{tabular}{@{}lllllllllllll@{}}
\toprule
\multicolumn{13}{c}{\textbf{Ranking Method}} \\ \midrule
\textbf{Instance Name} & \textbf{Ridge} & \textbf{Lasso} & \textbf{SVR} & \textbf{KNN} & \textbf{RF} & \textbf{MLP} & \textbf{Stacking} & \textbf{Voting} & \textbf{RBA} & \textbf{MW} & \textbf{GCGOS} & \textbf{Goodness} \\ \midrule
cost266-UUE.mps & 0.126 & 0.411 & 0.642 & 0.265 & \textbf{0.000} & 0.206 & 0.434 & 0.378 & 0.363 & 0.642 & 0.363 & 0.642 \\
dfn-bwin-DBE.mps & 0.059 & 0.483 & 0.035 & 0.572 & 0.057 & 0.154 & 0.405 & \textbf{0.020} & 0.483 & 0.510 & 0.483 & 0.410 \\
germany50-UUM.mps & 0.704 & 0.311 & \textbf{0.071} & 0.226 & 0.147 & 0.037 & 0.239 & \textbf{0.071} & \textbf{0.071} & 0.704 & \textbf{0.071} & 0.704 \\
ta1-UUM.mps & 0.083 & 0.329 & 0.342 & 0.139 & 0.156 & 0.284 & 0.075 & 0.146 & \textbf{0.000} & 0.438 & 0.017 & 0.438 \\
ta2-UUE.mps & \textbf{0.132} & \textbf{0.132} & \textbf{0.132} & 0.736 & \textbf{0.132} & \textbf{0.132} & 0.631 & \textbf{0.132} & 0.416 & 0.808 & 0.416 & 0.808 \\ \midrule
g200x740.mps & 0.197 & \textbf{0.000} & 0.149 & 0.128 & 0.088 & 0.145 & 0.108 & 0.138 & \textbf{0.000} & 0.435 & \textbf{0.000} & 0.250 \\
h50x2450.mps & \textbf{0.000} & \textbf{0.000} & \textbf{0.000} & 0.198 & \textbf{0.000} & 0.182 & 0.173 & 0.182 & 0.372 & 0.824 & 0.372 & 0.824 \\
h80x6320d.mps & \textbf{0.000} & \textbf{0.000} & \textbf{0.000} & \textbf{0.000} & \textbf{0.000} & \textbf{0.000} & \textbf{0.000} & \textbf{0.000} & \textbf{0.000} & 0.829 & \textbf{0.000} & 0.829 \\
k16x240b.mps & 0.319 & \textbf{0.000} & 0.188 & 0.454 & 0.525 & 0.136 & 0.147 & 0.176 & 0.035 & 0.826 & 0.035 & 0.858 \\ \midrule
snp-02-004-104.mps & 0.512 & 0.148 & 0.138 & \textbf{0.091} & 0.108 & 0.138 & 0.172 & \textbf{0.091} & 0.148 & 0.110 & 0.148 & \textbf{0.091} \\
snp-04-052-052.mps & 0.592 & 0.530 & 0.348 & 0.310 & \textbf{0.036} & 0.425 & 0.179 & 0.179 & 0.530 & 0.220 & 0.530 & 0.220 \\
snp-06-004-052.mps & 0.205 & 0.147 & 0.088 & 0.068 & 0.192 & 0.191 & 0.157 & 0.157 & 0.494 & \textbf{0.000} & 0.494 & \textbf{0.000} \\
snp-10-004-052.mps & 0.244 & 0.196 & \textbf{0.000} & 0.018 & 0.090 & \textbf{0.000} & \textbf{0.000} & \textbf{0.000} & 0.499 & 0.133 & 0.499 & 0.133 \\
snp-10-052-052.mps & 0.072 & 0.242 & \textbf{0.000} & 0.210 & 0.032 & \textbf{0.000} & \textbf{0.000} & \textbf{0.000} & 0.549 & 0.032 & 0.549 & 0.032 \\ \midrule
Average & 0.232 & 0.209 & 0.152 & 0.244 & 0.112 & 0.145 & 0.194 & 0.119 & 0.283 & 0.465 & 0.284 & 0.446 \\  \bottomrule
\end{tabular}
\end{table*}
\end{landscape}

\subsubsection{The Best Ranking Method}
We compared all ranking methods against one another using the Friedman test statistic which can be calculated as $F_f = \frac{12n}{k(k+1)} [\sum_{j}R^2_{j} - \frac{k(k+1)^2}{4}]$, where $n$ is the number of test instances, $k$ is the number of comparison algorithms, $R_{j}$ is the average rank of algorithm $j$. Using the CONTROLTEST package as referenced in \citep{Derrac2011} showed a statistical significance between the different ranking methods ($p = 0.0133$). To reduce the chance of making a Type-1 error when comparing multiple algorithms, we carried out pairwise comparisons between the best ML method (Voting) and the best heuristic methods (RBA and Goodness) using the Aligned Friedman Ranks. These heuristic methods were chosen as they produce similar prediction results to the GCGOS and MW methods respectively, albeit slightly better. These pairwise comparisons give associated p-values of $0.033$ and $3.81e^{-04}$ respectively, showing statistical significance at the $\alpha = 0.05$ level and indicating the superiority of a ML ranking method at predicting decomposition qualities. Amongst all ranking functions the Voting method was able to provide equal best predictions for 6 out of the 14 instances and the outright best prediction for an additional instance. The Voting method was also shown to never predict the worst decomposition when compared to the heuristic ranking techniques. It should be noted that the performances of the different ML models were not significantly different from one another, indicating that further parameter tuning is unlikely to yield significantly better results. In addition, one of the key undertakings of this paper was to determine if a ML approach in general can more accurately predict decomposition quality than traditional heuristic measures, which has been demonstrated in our findings. For these reasons we have not carried out any extensive hyper-parameter tuning, although future work may look at carrying out hyper-parameter tuning on a selected ML model in an attempt to further improve ranking performance.

\subsubsection{Improvements over a Multi-Objective Approach}
\begin{figure}[ht]
\centering
  \subfloat[h50x2450.mps]{
   \includegraphics[height=0.29\textwidth]{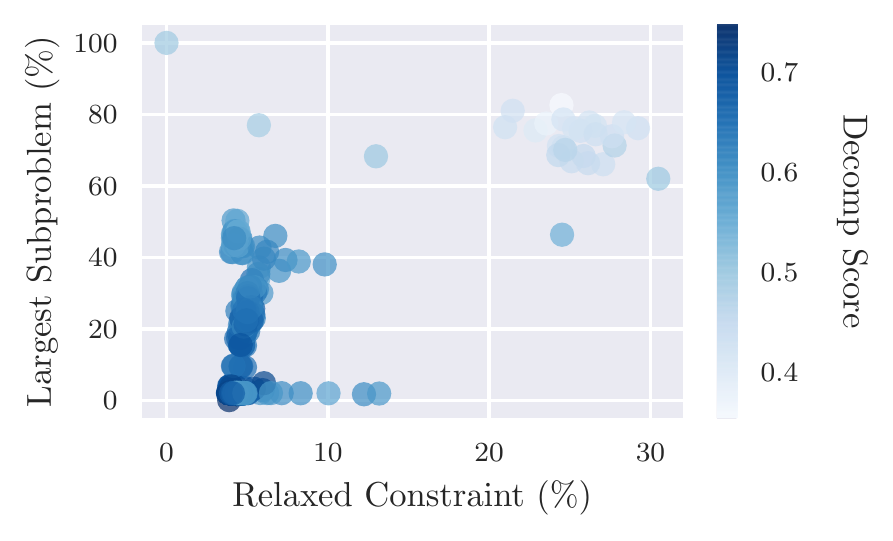}
   }
 \hfill 	
  \subfloat[snp-10-004-052.mps]{
   \includegraphics[height=0.29\textwidth]{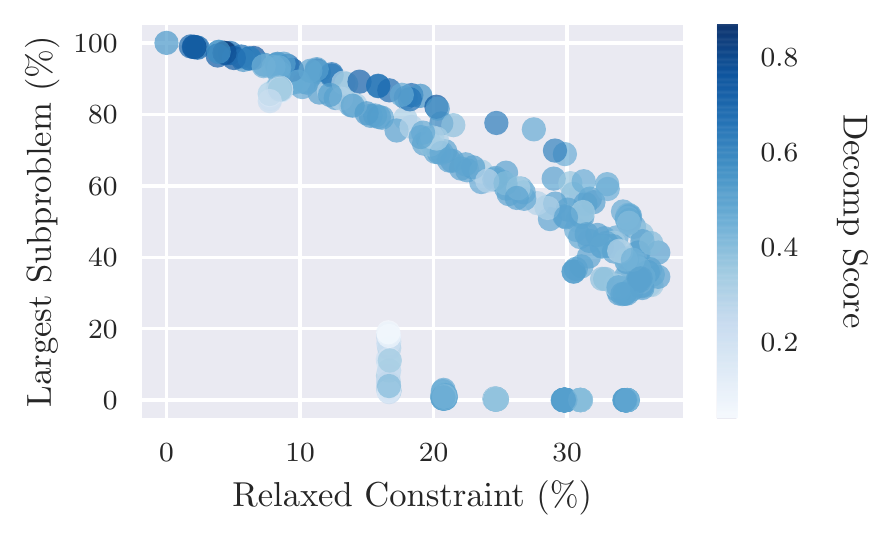}
   }
\caption{Multi-objective Correlation: Shown for all decompositions in both h50x2450.mps and snp-10-004-052.mps instances are the correlations between Relaxed Constraint (\%) and Largest Subproblem (\%), where Largest Subproblem (\%) is calculated as the percentage of all MIP variables contained within the largest subproblem. As seen in h50x2450.mps, minimising both the Relaxed Constraint (\%) and the Largest Subproblem (\%) does not correlate well decomposition quality. In contrast, the snp-10-004-052.mps instance shows good correlation between minimising both the given objectives and decomposition quality, although there are many good decompositions both in and outside of the Pareto optimal solution front.}
\label{fig:MultiobjectiveCorrelation}
\end{figure}

\begin{figure}[ht]
\centering
  \subfloat[h50x2450.mps]{
   \includegraphics[width=0.48\textwidth]{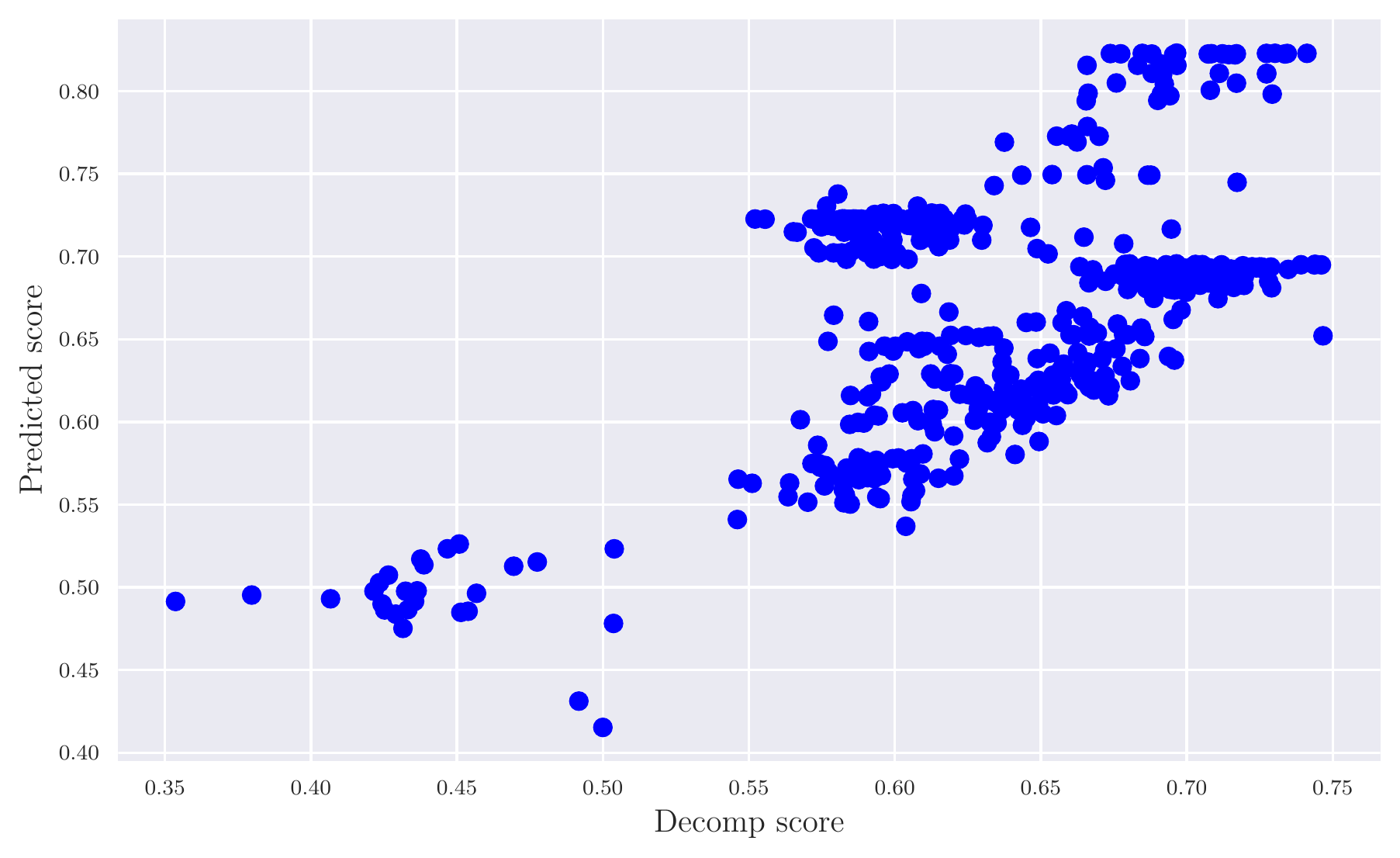}
}
 \hfill 	
  \subfloat[snp-10-004-052.mps]{
   \includegraphics[width=0.48\textwidth]{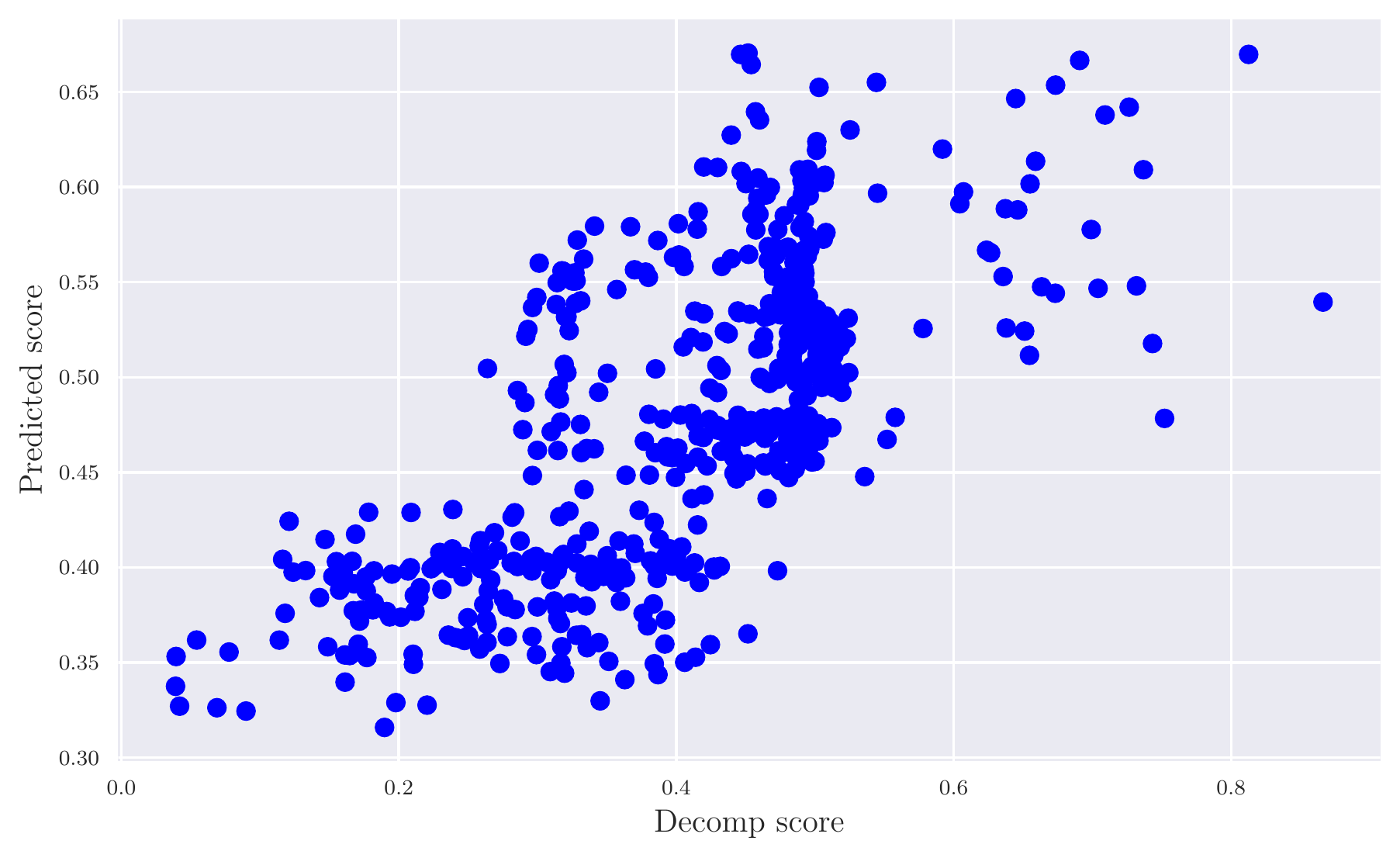}
}
\caption{Machine Learning Correlation: Shown for all decompositions in both h50x2450.mps and snp-10-004-052.mps instances are the correlations between the predicted decomposition score and the actual decomposition score for the trained Voting ensemble function. As seen in both h50x2450.mps and snp-10-004-052.mps instances, the ML ranking function shows good correlation between predicted decomposition scores and actual decomposition scores.}
\label{fig:PredictionCorrelation}
\end{figure}

Whilst a population based metaheuristic such as the NSGA-II algorithm is able to produce a large pool of solutions in a reasonable amount of time, the current objectives used to evolve the population might not always correlate well with decomposition quality. In addition, as demonstrated in \citep{Weiner2020}, even amongst the final Pareto optimal solutions there can be significant variability in bound quality and run time, requiring the practitioner to then test a potentially large number of decompositions. Shown in Figure~\ref{fig:MultiobjectiveCorrelation} and Figure~\ref{fig:PredictionCorrelation} are two examples of where a ranking function could be applied in conjunction with an evolutionary algorithm to identify high quality decompositions, without having to evaluate all Pareto optimal solutions found by the NSGA-II algorithm. In Figure~\ref{fig:MultiobjectiveCorrelation}(a) there appears to be no correlation between decomposition quality and the two objectives used in the NSGA-II algorithm, highlighting when considering only these two features results in poor decompositions. In Figure~\ref{fig:MultiobjectiveCorrelation}(b), whilst there is a correlation between decomposition quality and the two objectives used in the NSGA-II algorithm, there are many high quality decompositions which lie outside of the Pareto front, and would potentially be ignored. Using a ML approach however, shows a high correlation between predicted decomposition scores and actual decomposition scores as shown in Figure~\ref{fig:PredictionCorrelation}. The ML ranking function can therefore be suitable in selecting which decompositions to run through a full LR framework, or can be used as a search guide instead of the two objectives used in the NSGA-II algorithm.


\subsection{Considering Relaxed Constraint Features Only} \label{Subsection: Relaxed Constraint Features Only}
\begin{figure*}[!htb]
\centering
\includegraphics[width=0.7\linewidth]{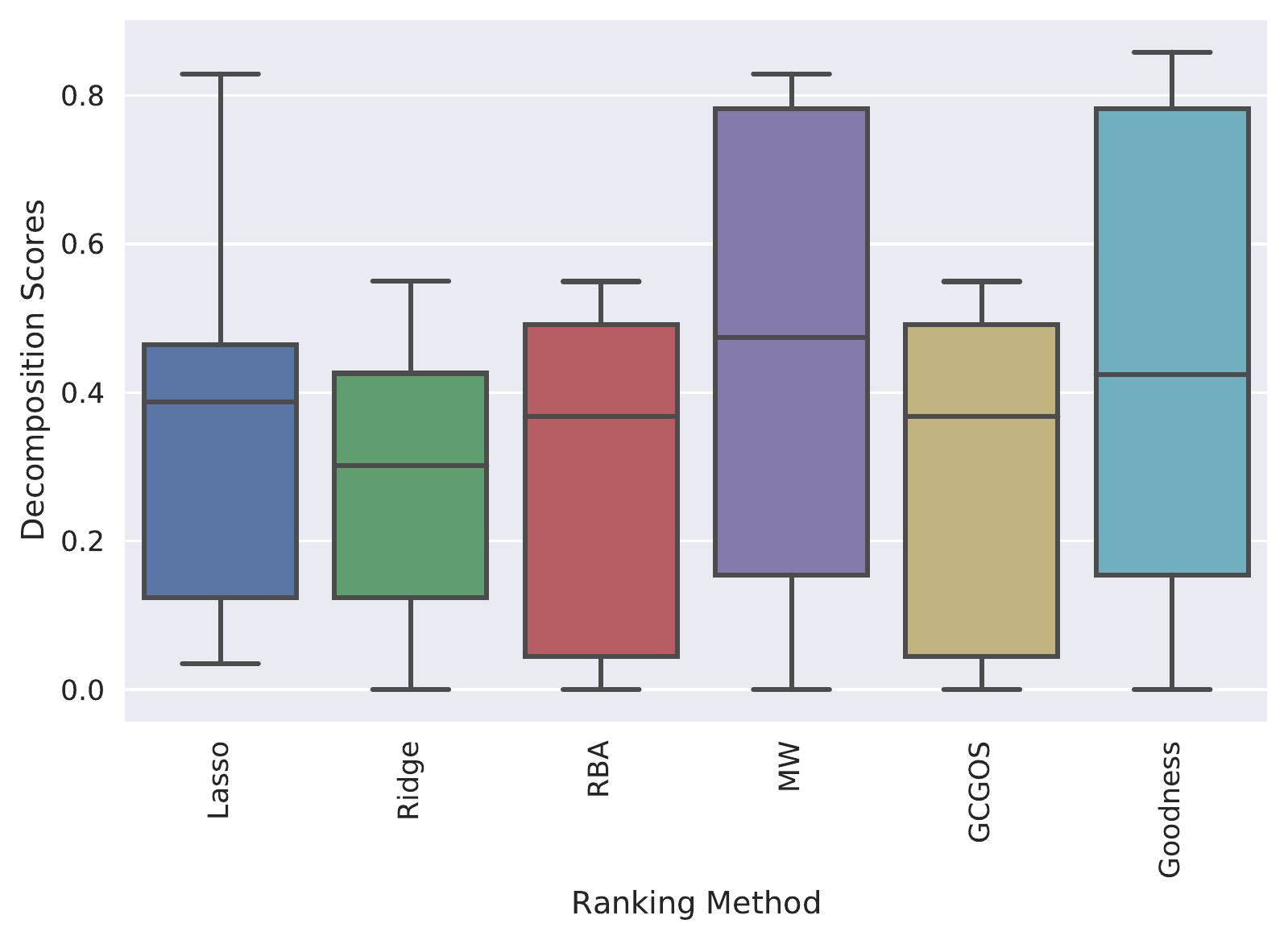}\vspace{-2mm}
\caption{Relaxed Constraint Features and Prediction Quality. Shown for each ranking method is a boxplot of the scores for the best decompositions selected by each ranking method across the 14 test instances in the Network dataset. The two Linear Regression models using Lasso and Ridge regularisation were trained using only the selected relaxed constraint features. These best selected decompositions by the ML methods using only the relaxed constraint features are still competitive with the decompositions selected by the heuristic ranking methods. The decomposition scores range from 0 (the best decomposition in the test set) to 1 (the worst decomposition in the test set).}
\label{fig:network_relaxed_constraint_predictions}
\end{figure*}

Currently, in the greedy random selection procedure used to generate the decompositions in this paper, constraints are selected for relaxation in a probabilistic manner using only non-zero proportions. Through manual feature selection, we considered three additional features to create a better constraint ranking function, including the average binary + integer proportions of relaxed constraints, the average RHS value of relaxed constraints and the average sum of objective coefficients of the variables associated with relaxed constraints. These features were chosen as they are easily calculable and can easily be incorporated in a greedy constraint selection process for future work. We tested these features and their effect on prediction quality by training two Linear Regression models with Lasso and Ridge regularisation respectively. A Linear Regression model was chosen as such a model is transparent, containing feature coefficients which are easily extractable to use in a future constraint ranking function. As shown in Figure~\ref{fig:network_relaxed_constraint_predictions}, for the Network dataset, a Linear Regression model using only these four features is often able to produce prediction qualities better than or equal to the heuristic techniques, showing promise that such a model has good predictive capability. Shown in Table~\ref{tab:rc_ridge_feature_coefficients} are the coefficient values found in the Linear Regression model using Ridge regularisation.

\begin{table}[!htb]
\centering
\caption{Relaxed Constraint Feature Coefficients. Shown for the Linear Regression model with Ridge Regularisation are the coefficients of the features in the model when trained on the Network dataset. These coefficients can be used for future constraint relaxation selection procedures, as this model has demonstrated good predictive capabilities.}
\label{tab:rc_ridge_feature_coefficients}
\begin{tabular}{@{}lr@{}}
\toprule
\textbf{Feature} & \textbf{Coefficient} \\ \midrule
Average\_Relaxed\_Constraint\_Statistics\_Non\_zero\_props & -10.238 \\
Average\_Relaxed\_Constraint\_Statistics\_Sum\_obj & -0.180 \\
Average\_Relaxed\_Constraint\_Statistics\_RHS\_vals & 0.253 \\
Average\_Relaxed\_Constraint\_Statistics\_Bin\_Int\_props & 4.223 \\ \bottomrule
\end{tabular}
\end{table}



\subsection{Testing on Randomly Selected Instances}

\begin{figure*}[!th]
\centering
\includegraphics[width=0.7\linewidth]{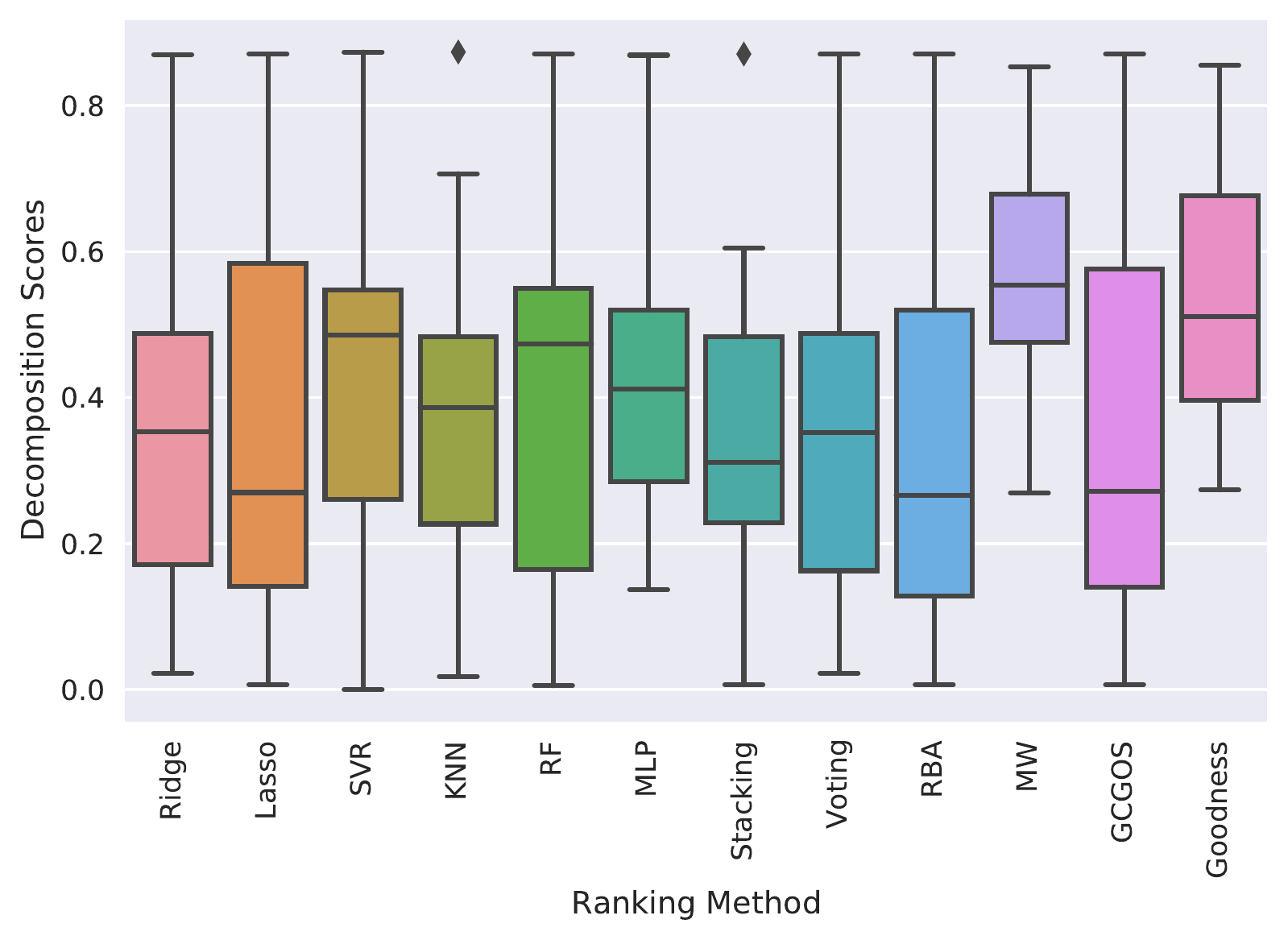}\vspace{-2mm}
\caption{Random MIPLIB Prediction Benchmarking. Shown for each ranking method is a boxplot of the scores for the best decompositions selected by each ranking method across the 10 random MIPLIB instances. The ML models presented were trained on all instances from the Network dataset. The decomposition scores range from 0 (the best decomposition in the test set) to 1 (the worst decomposition in the test set).}
\label{fig:random_miplib_boxplot}
\end{figure*}

\begin{figure}[!th]
\centering
\includegraphics[width=0.9\textwidth]{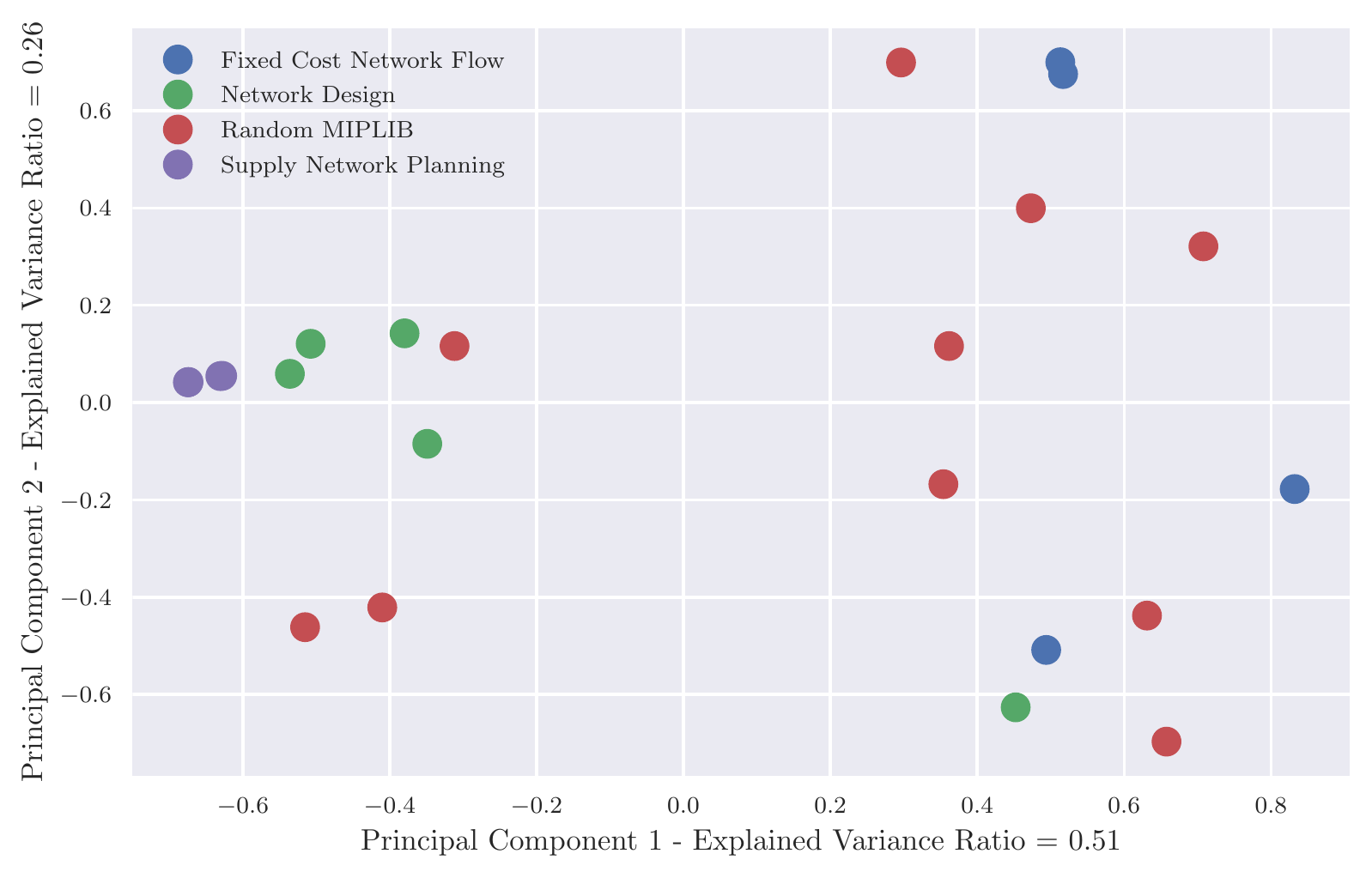}\vspace{-2mm}
\caption{PCA Analysis of all problem types including Network Design, Fixed Cost Network Flow, Supply Networking and Random MIPLIB instances. Shown in this figure are the first two principle components, comprising over 77\% of the explained variance. As can be seen, the Random MIPLIB instances appear to show no clustering patterns in comparison to the Network subproblem types such as Supply Network Planning and Network Design types.}
\label{fig:PCA All Problems}
\end{figure}

A final experiment was carried out to see how a ML model trained on all instances from the Network dataset would perform on randomly selected instances from the MIPLIB 2017 library. The prediction results for the different ranking methods considered are presented in Figure~\ref{fig:random_miplib_boxplot}. As can be seen, the ML based ranking functions are competitive with the heuristic based techniques without being state-of-the-art. Using the Friedman test, no statistical significance was detected ($p = 0.419$) amongst the different ranking methods. These findings appear to be consistent with those in \citep{Basso2020}, indicating that a ML based ranking function is only useful when test instances are somewhat similar to the training instances. A PCA plot using the instance features described in Table~A7 is shown in Figure~\ref{fig:PCA All Problems}. Unlike previous problem types tested, the randomly selected MIPLIB instances do not seem to form any significant clusters, showing that at least in this feature space the instances are quite dissimilar. 

\subsection{Discussion}

In light of our findings, we present a brief discussion of how our work could be leveraged in future automatic decomposition frameworks. As demonstrated in our findings, there is a clear link between instance similarity and the predictive capabilities of a ML based ranking function for predicting decomposition qualities. As such, we believe that future works relying on decomposition quality predictions should have an initial instance classification procedure, selecting a pre-trained ML model that might be most suitable for the instance to be solved. Whilst outside the scope of this paper, this classification step requires significantly more data and computational work for new problem types not explored in this paper, in order to create a suitable set of decompositions for training purposes. In addition to classifying instances based on problem type, classification on other metrics might also be appropriate, such as the proportions of different constraint types e.g., Set Partitioning, Binpacking, Set Covering, Knapsack etc. Unlike in previous studies \citep{Basso2020, Weiner2020}, where constraints are relaxed in a greedy manner using only their non-zero proportions, we propose that other constraint features be considered such as those shown in Table~\ref{tab:rc_ridge_feature_coefficients}, as we have demonstrated the importance these features have on decomposition quality. Generating a valid decomposition via constraint relaxation can be carried out significantly fast using a hypergraph partitioning method as discussed in \citep{Bergner2015a, Weiner2020}, requiring $O(nz)$ time complexity, where $nz$ is the number of non-zeroes in the constraint matrix. Such a phenomenon means that a large population of decompositions can be ranked using a ML function relatively quickly, in order to determine which decompositions should be selected for solving, or potentially as a candidate for further improvements using a local search operator.

\clearpage

\section{Conclusions}\label{Section:Conclusion}
This paper explored how a Machine Learning (ML) approach can be used effectively as a decomposition ranking function for Lagrangian Relaxation of Mixed Integer Programs. Through benchmark comparisons with previously published hand-crafted heuristic measures, it was demonstrated that a ML approach was able to provide state-of-the-art results in predicting decomposition qualities from amongst a large number of candidate decompositions. This paper produced a rich dataset of decompositions relating to both network type instances and randomly selected instances from the MIPLIB 2017 library which is freely available for other researchers to access \footnote{The dataset and code used in this paper will become freely accessible upon the acceptance of this paper.}. This paper also explored how instance similarity plays a critical role in ML prediction qualities, suggesting that future work in instance classification and the exploration of additional instance features could provide a promising research direction. Finally, a new constraint ranking function was also provided, which has shown promising prediction capabilities and can be used to rank constraints for relaxation a-priori to any heuristic selection algorithm. Future work may involve embedding the ML ranking function described in this paper within a heuristic based search technique, in order to find high quality solutions without directly evaluating a large population of candidate decompositions. Whilst better results could potentially be found through hyperparameter tuning of models, by primarily using only default model settings we were still able to demonstrate that a ML approach can significantly outperform current benchmark heuristic methods. Future work may look at testing both more ML models as well as additional hyperparameter tuning.



\section*{Acknowledgements}
This research was supported by an ARC (Australian Research Council) Discovery Grant (DP180101170)


\bibliography{Machine_Learning.bib}
\end{document}


\begin{table}[!t]
\setlength\tabcolsep{30pt}
\def\arraystretch{2}
\centering
\caption{Full Features List: Shown is the full list of features used to train and test the Machine Learning models}
\label{tab:Full_Features_List}
\begin{tabular}{@{}l@{}}
\toprule
\textbf{Relaxed Constraint Statistics} \\ \midrule
Min, Max, Ave, Stddev of Constraint RHS Value \\
Min, Max, Ave, Stddev of Sum of Obj Coefficients of Variables in Constraint \\
Min, Max, Ave, Stddev of Prop of Instance Bin in Constraint \\
Min, Max, Ave, Stddev of Prop of Instance Int in Constraint \\
Min, Max, Ave, Stddev of Prop of Instance Bin + Int in Constraint \\
Min, Max, Ave, Stddev of Non Zero Props in Constraint \\ 
Prop of Constraints Relaxed \\
Prop of Constraints Relaxed which are Equality \\ \midrule
\textbf{Subproblem Statisitics} \\ \midrule
Min, Max, Ave, Stddev of Prop of Instance Bin in Subproblems \\
Min, Max, Ave, Stddev of Prop of Instance Int in Subproblems \\
Min, Max, Ave, Stddev of Prop of Instance Bin + Int in Subproblems \\
Min, Max, Ave, Stddev of Subproblem Densities \\
Min, Max, Ave, Stddev of Constraint Prop in Subproblems \\
Min, Max, Ave, Stddev Equality Prop in Subproblems \\
Min, Max, Ave, Stddev Variable Prop in Subproblems \\
Min, Max, Ave, Stddev Non Zero Props in Subproblems \\
Ave, Stddev of Sum of Objective Coefficients for Variables in Subproblems \\
Ave, Stddev of Range of Objective Coefficients for Variables in Subproblems \\
Ave, Stddev of Average RHS Value in Subproblems \\
Ave, Stddev of RHS Range in Subproblems \\
Ave, Stddev of Subproblem Shapes \\ \bottomrule
\end{tabular}
\end{table}

\begin{table}[!t]
\centering
\caption{Decomposition Counts: Shown for instance in the dataset are the total number decompositions created for both training and testing purposes.}
\label{tab:Decomposition Counts}
\begin{tabular}{@{}llr@{}}
\toprule
Problem Type            & Instance Name           & Number of Decompositions \\ \midrule
Network Design          & cost266-UUE.mps         & 1743                     \\
Network Design          & dfn-bwin-DBE.mps        & 1996                     \\
Network Design          & germany50-UUM.mps       & 1780                     \\
Network Design          & ta1-UUM.mps             & 1848                     \\
Network Design          & ta2-UUE.mps             & 1673                     \\ \midrule
Fixed Cost Network Flow & g200x740.mps            & 2268                     \\
Fixed Cost Network Flow & h50x2450.mps            & 504                      \\
Fixed Cost Network Flow & h80x6320d.mps           & 469                      \\
Fixed Cost Network Flow & k16x240b.mps            & 640                      \\ \midrule
Supply Network Planning & snp-02-004-104.mps      & 2485                     \\
Supply Network Planning & snp-04-052-052.mps      & 2462                     \\
Supply Network Planning & snp-06-004-052.mps      & 1538                     \\
Supply Network Planning & snp-10-004-052.mps      & 575                      \\
Supply Network Planning & snp-10-052-052.mps      & 582                      \\ \midrule
Random MIPLIB           & blp-ic98.mps            & 702                      \\
Random MIPLIB           & dws008-01.mps           & 2314                     \\
Random MIPLIB           & 30n20b8.mps             & 1669                     \\
Random MIPLIB           & air03.mps               & 2616                     \\
Random MIPLIB           & traininstance2.mps      & 2090                     \\
Random MIPLIB           & neos-4387871-tavua.mps  & 2535                     \\
Random MIPLIB           & neos-4338804-snowy.mps  & 1934                     \\
Random MIPLIB           & air05.mps               & 2607                     \\
Random MIPLIB           & neos-4954672-berkel.mps & 1841                     \\
Random MIPLIB           & splice1k1.mps           & 1318                     \\ \bottomrule
\end{tabular}
\end{table}

\begin{table}[!t]
\setlength\tabcolsep{2pt}
\centering
\caption{Raw Bound Results: Shown for each instance are the bounds of the min and max bounds found amongst all decompositions tested. As all instances solved are minimization type problems, smaller LR bounds are of worse quality than larger LR bounds. Shown for each instance is also the LP bound as well as the best known primal solution as reported on MIPLIB2017. It should be noted that as subproblems were solved to within 1\% of optimality, as the primal solutions to the subproblems were used in bound calculations when solutions were within this optimality tolerance, it is possible for the Max LR bound to be slightly higher (\textless 1\%) than the best known primal solution.}
\label{tab:Raw Bound Results}
\begin{tabular}{@{}lrrrr@{}}
                        & \multicolumn{4}{c}{\textbf{Bound Comparisons}}                \\ \midrule
\textbf{Instance} & \textbf{Min LR Bound} & \textbf{Max LR Bound} & \textbf{LP Bound} & \textbf{Best Known Primal Solution} \\ \midrule
cost266-UUE.mps         & 20161500.00   & 24150800.00   & 20161500.00   & 25148940.56   \\
dfn-bwin-DBE.mps        & 17890.90      & 50921.20      & 17890.90      & 73623.79      \\
germany50-UUM.mps       & 597932.00     & 618552.00     & 597932.00     & 628490.00     \\
ta1-UUM.mps             & 3693670.00    & 7237730.00    & 3693670.00    & 7518328.20    \\
ta2-UUE.mps             & 36964000.00   & 37886500.00   & 36964000.00   & 37871728.59   \\ \midrule
g200x740.mps            & 34077.50      & 44356.00      & 34077.50      & 44316.00      \\
h50x2450.mps            & 11147.70      & 32972.50      & 11147.70      & 32906.88      \\
h80x6320d.mps           & 5325.16       & 6382.10       & 5325.16       & 6382.10       \\
k16x240b.mps            & 3320.77       & 11331.70      & 3320.77       & 11393.00      \\ \midrule
snp-02-004-104.mps      & 548045000.00  & 586972000.00  & 548045000.00  & 586803238.66  \\
snp-04-052-052.mps      & 728302000.00  & 857416000.00  & 728196000.00  & 885202237.19  \\
snp-06-004-052.mps      & 1787140000.00 & 1875140000.00 & 1787140000.00 & 1869531919.90 \\
snp-10-004-052.mps      & 5842870000.00 & 5914170000.00 & 5842630000.00 & 5906642865.78 \\
snp-10-052-052.mps      & 5843240000.00 & 5944030000.00 & 5842630000.00 & 6364531568.74 \\ \midrule
blp-ic98.mps            & 4331.17       & 4515.02       & 4331.17       & 4491.45       \\
dws008-01.mps           & 584.50        & 25532.30      & 584.50        & 37412.60      \\
30n20b8.mps             & 1.57          & 302.00        & 1.57          & 302.00        \\
air03.mps               & 338864.00     & 342760.00     & 338864.00     & 340160.00     \\
traininstance2.mps      & 0.00          & 17760.00      & 0.00          & 71820.00      \\
neos-4387871-tavua.mps  & 10.14         & 27.85         & 10.14         & 33.38         \\
neos-4338804-snowy.mps  & 1447.00       & 1473.00       & 1447.00       & 1471.00       \\
air05.mps               & 25877.60      & 26497.00      & 25877.60      & 26374.00      \\
neos-4954672-berkel.mps & 1150230.00    & 2345130.00    & 1150230.00    & 2612710.00    \\
splice1k1.mps           & -1646.78      & -798.09       & -1646.78      & -394.00       \\ \bottomrule
\end{tabular}
\end{table}

\begin{landscape}
\begin{figure}[!hbt]
\centering
\includegraphics[width=\linewidth]{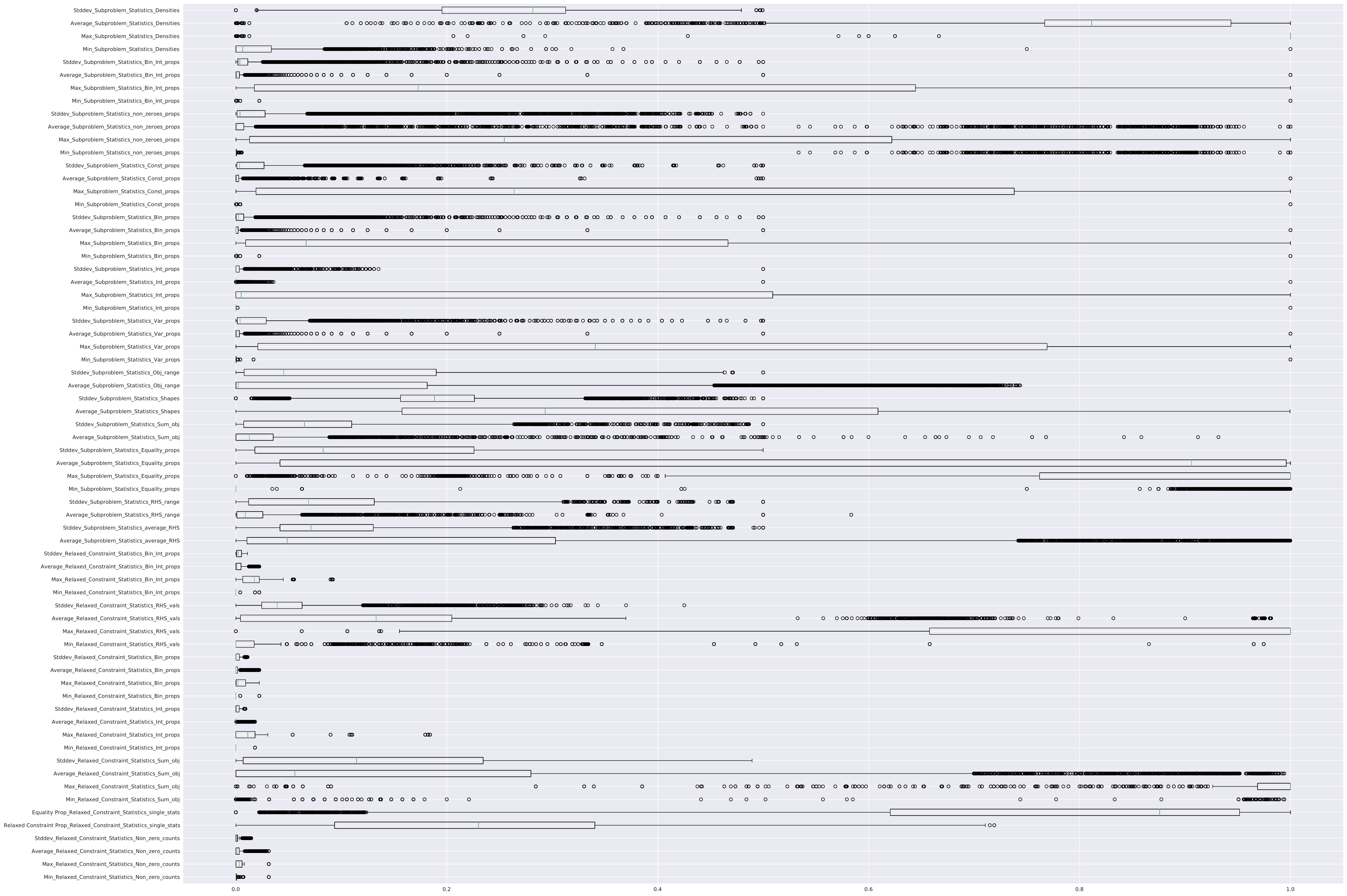}\vspace{-2mm}
\caption{Network Dataset Feature Spread. Shown are the spread of feature data for all decompositions contained in the Network dataset, containing the network problem types - Network Design, Fixed Cost Network Flow and Supply Network Planning. In total there are 20563 decompositions in the Network dataset obtained from 14 instances.}
\label{fig:Network Feature Spread}
\end{figure}
\end{landscape}

\begin{landscape}
\begin{figure}[!hbt]
\centering
\includegraphics[width=\linewidth]{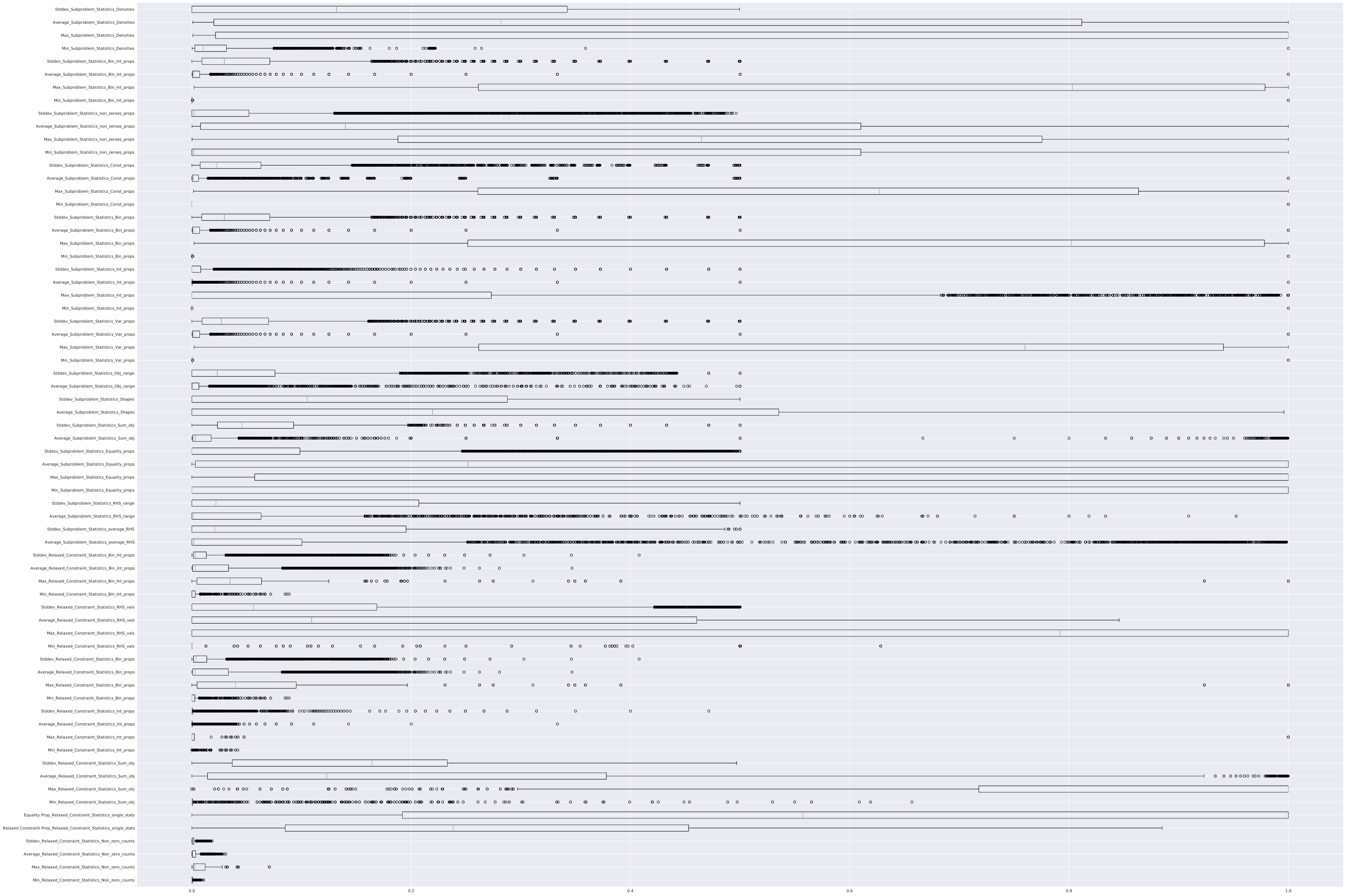}\vspace{-2mm}
\caption{Random Dataset Decomposition Feature Spread. Shown are the spread of feature data for all decompositions contained in the Random dataset. In total there are 19626 decompositions in the Random dataset obtained from 10 instances.}
\label{fig:Random Feature Spread}
\end{figure}
\end{landscape}

\begin{table}[]
\centering
\caption{Network Design Prediction Results: Shown for each ML ranking method and test instance (Model\_Instance) is the score for the best decomposition selected for each instance from the Network Design problem type, with 0 indicating the best decomposition from the test instance was selected and 1 indicating the worst decomposition from the test instance was selected. Each column represents the problem type the ML model was trained on. For the network design problem type, the model was trained on all instances except for the test instance. For the other problem types, the model was trained on all instances. To validate the significance the problem type the model is trained on has on performance, the z-scores and associated p-values shown are calculated via pairwise comparisons between the control method (models trained on the same problem type) and comparison methods (models trained on different problem types) using the Aligned Friedman ranks. Finally, the average score for the best decompositions selected by each model when trained on the different problem types is also presented.}
\label{tab:Intersubuproblem-ND}
\begin{tabular}{@{}lrrr@{}}
\toprule
Model\_Instance & network\_design & fixed\_cost\_network\_flow & supply\_network\_planning \\ \midrule
Lasso\_cost266-UUE.mps & 0.000 & 0.163 & 0.628 \\
Lasso\_dfn-bwin-DBE.mps & 0.415 & 0.483 & 0.410 \\
Lasso\_germany50-UUM.mps & 0.071 & 0.244 & 0.239 \\
Lasso\_ta1-UUM.mps & 0.017 & 0.206 & 0.054 \\
Lasso\_ta2-UUE.mps & 0.594 & 0.132 & 0.800 \\
Ridge\_cost266-UUE.mps & 0.231 & 0.628 & 0.040 \\
Ridge\_dfn-bwin-DBE.mps & 0.415 & 0.483 & 0.415 \\
Ridge\_germany50-UUM.mps & 0.071 & 0.704 & 0.071 \\
Ridge\_ta1-UUM.mps & 0.167 & 0.437 & 0.000 \\
Ridge\_ta2-UUE.mps & 0.132 & 0.800 & 0.391 \\
SVR\_cost266-UUE.mps & 0.000 & 0.206 & 0.642 \\
SVR\_dfn-bwin-DBE.mps & 0.109 & 0.273 & 0.000 \\
SVR\_germany50-UUM.mps & 0.077 & 0.192 & 0.704 \\
SVR\_ta1-UUM.mps & 0.306 & 0.358 & 0.438 \\
SVR\_ta2-UUE.mps & 0.087 & 0.132 & 0.391 \\
KNN\_cost266-UUE.mps & 0.265 & 0.008 & 0.364 \\
KNN\_dfn-bwin-DBE.mps & 0.440 & 0.096 & 0.493 \\
KNN\_germany50-UUM.mps & 0.071 & 0.259 & 0.698 \\
KNN\_ta1-UUM.mps & 0.211 & 0.270 & 0.438 \\
KNN\_ta2-UUE.mps & 0.736 & 0.132 & 0.449 \\
RF\_cost266-UUE.mps & 0.185 & 0.466 & 0.642 \\
RF\_dfn-bwin-DBE.mps & 0.619 & 0.285 & 0.562 \\
RF\_germany50-UUM.mps & 0.130 & 0.199 & 0.244 \\
RF\_ta1-UUM.mps & 0.154 & 0.084 & 0.395 \\
RF\_ta2-UUE.mps & 0.381 & 0.800 & 0.800 \\
MLP\_cost266-UUE.mps & 0.100 & 0.040 & 0.078 \\
MLP\_dfn-bwin-DBE.mps & 0.037 & 0.571 & 0.101 \\
MLP\_germany50-UUM.mps & 0.037 & 0.077 & 0.071 \\
MLP\_ta1-UUM.mps & 0.182 & 0.182 & 0.017 \\
MLP\_ta2-UUE.mps & 0.632 & 0.407 & 0.393 \\
Stacking\_cost266-UUE.mps & 0.158 & 0.040 & 0.078 \\
Stacking\_dfn-bwin-DBE.mps & 0.641 & 0.571 & 0.402 \\
Stacking\_germany50-UUM.mps & 0.071 & 0.077 & 0.071 \\
Stacking\_ta1-UUM.mps & 0.167 & 0.182 & 0.017 \\
Stacking\_ta2-UUE.mps & 0.132 & 0.132 & 0.304 \\
Voting\_cost266-UUE.mps & 0.116 & 0.475 & 0.078 \\
Voting\_dfn-bwin-DBE.mps & 0.179 & 0.483 & 0.421 \\
Voting\_germany50-UUM.mps & 0.071 & 0.366 & 0.071 \\
Voting\_ta1-UUM.mps & 0.167 & 0.424 & 0.000 \\
Voting\_ta2-UUE.mps & 0.632 & 0.449 & 0.391 \\ \midrule
Average & \textbf{0.230} & 0.313 & 0.320 \\ \midrule
z-score &  & -3.531 & -3.675 \\
p-value &  & 2.07E-04 & 1.19E-04 \\ \bottomrule
\end{tabular}
\end{table}

\begin{table}[]
\centering
\caption{Fixed Cost Network Flow Prediction Results: Shown for each ML ranking method and test instance (Model\_Instance) is the score for the best decomposition selected for each instance from the Fixed Cost Network Flow problem type, with 0 indicating the best decomposition from the test instance was selected and 1 indicating the worst decomposition from the test instance was selected. Each column represents the problem type the ML model was trained on. For the Fixed Cost Network Flow, the model was trained on all instances except for the test instance. For the other problem types, the model was trained on all instances. To validate the significance the problem type the model is trained on has on performance, the z-scores and associated p-values shown are calculated via pairwise comparisons between the control method (models trained on the same problem type) and comparison methods (models trained on different problem types) using the Aligned Friedman ranks. Finally, the average score for the best decompositions selected by each model when trained on the different problem types is also presented.}
\label{tab:Intersubuproblem-FCNF}
\begin{tabular}{@{}lrrr@{}}
\toprule
Model\_Instance & fixed\_cost\_network\_flow & network\_design & supply\_network\_planning \\ \midrule
Lasso\_g200x740.mps & 0.000 & 0.000 & 0.000 \\
Lasso\_h50x2450.mps & 0.000 & 0.000 & 0.372 \\
Lasso\_h80x6320d.mps & 0.000 & 0.000 & 0.000 \\
Lasso\_k16x240b.mps & 0.152 & 0.000 & 0.000 \\
Ridge\_g200x740.mps & 0.113 & 0.077 & 0.000 \\
Ridge\_h50x2450.mps & 0.066 & 0.000 & 0.000 \\
Ridge\_h80x6320d.mps & 0.073 & 0.000 & 0.000 \\
Ridge\_k16x240b.mps & 0.107 & 0.759 & 0.035 \\
SVR\_g200x740.mps & 0.137 & 0.139 & 0.000 \\
SVR\_h50x2450.mps & 0.000 & 0.000 & 0.173 \\
SVR\_h80x6320d.mps & 0.000 & 0.000 & 0.000 \\
SVR\_k16x240b.mps & 0.107 & 0.497 & 0.107 \\
KNN\_g200x740.mps & 0.108 & 0.108 & 0.176 \\
KNN\_h50x2450.mps & 0.372 & 0.000 & 0.173 \\
KNN\_h80x6320d.mps & 0.000 & 0.073 & 0.074 \\
KNN\_k16x240b.mps & 0.035 & 0.233 & 0.175 \\
RF\_g200x740.mps & 0.118 & 0.077 & 0.541 \\
RF\_h50x2450.mps & 0.000 & 0.135 & 0.372 \\
RF\_h80x6320d.mps & 0.000 & 0.074 & 0.000 \\
RF\_k16x240b.mps & 0.789 & 0.187 & 0.136 \\
MLP\_g200x740.mps & 0.107 & 0.108 & 0.000 \\
MLP\_h50x2450.mps & 0.000 & 0.751 & 0.173 \\
MLP\_h80x6320d.mps & 0.074 & 0.848 & 0.000 \\
MLP\_k16x240b.mps & 0.107 & 0.000 & 0.035 \\
Stacking\_g200x740.mps & 0.118 & 0.077 & 0.000 \\
Stacking\_h50x2450.mps & 0.000 & 0.135 & 0.173 \\
Stacking\_h80x6320d.mps & 0.000 & 0.074 & 0.000 \\
Stacking\_k16x240b.mps & 0.000 & 0.793 & 0.035 \\
Voting\_g200x740.mps & 0.118 & 0.077 & 0.000 \\
Voting\_h50x2450.mps & 0.000 & 0.000 & 0.180 \\
Voting\_h80x6320d.mps & 0.000 & 0.000 & 0.000 \\
Voting\_k16x240b.mps & 0.107 & 0.455 & 0.035 \\ \midrule
Average & \textbf{0.088} & 0.177 & 0.093 \\ \midrule
z-score &  & -2.877 & 0.177 \\
p-value &  & 0.002 & 0.570 \\ \bottomrule
\end{tabular}
\end{table}

\begin{table}[]
\centering
\caption{Supply Network Planning Prediction Results: Shown for each ML ranking method and test instance (Model\_Instance) is the score for the best decomposition selected each instance from the Supply Network Planning problem type, with 0 indicating the best decomposition from the test instance was selected and 1 indicating the worst decomposition from the test instance was selected. Each column represents the problem type the ML model was trained on. For the Supply Network Planning problem type, the model was trained on all instances except for the test instance. For the other problem types, the model was trained on all instances. To validate the significance the problem type the model is trained on has on performance, the z-scores and associated p-valuesshown are calculated via pairwise comparisons between the control method (models trained on the same problem type) and comparison methods (models trained on different problem types) using the Aligned Friedman ranks. Finally, the average score for the best decompositions selected by each model when trained on the different problem types is also presented.}
\label{tab:Intersubuproblem-SNP}
\begin{tabular}{lrrr}
\hline
Model\_Instance & supply\_network\_planning & fixed\_cost\_network\_flow & network\_design \\ \midrule
Lasso\_snp-02-004-104.mps & 0.499 & 0.019 & 0.148 \\
Lasso\_snp-04-052-052.mps & 0.530 & 0.530 & 0.530 \\
Lasso\_snp-06-004-052.mps & 0.155 & 0.078 & 0.494 \\
Lasso\_snp-10-004-052.mps & 0.036 & 0.000 & 0.499 \\
Lasso\_snp-10-052-052.mps & 0.000 & 0.487 & 0.549 \\
Ridge\_snp-02-004-104.mps & 0.487 & 0.054 & 0.398 \\
Ridge\_snp-04-052-052.mps & 0.278 & 0.593 & 0.823 \\
Ridge\_snp-06-004-052.mps & 0.155 & 0.080 & 0.491 \\
Ridge\_snp-10-004-052.mps & 0.036 & 0.207 & 0.484 \\
Ridge\_snp-10-052-052.mps & 0.000 & 0.437 & 0.803 \\
SVR\_snp-02-004-104.mps & 0.138 & 0.012 & 0.148 \\
SVR\_snp-04-052-052.mps & 0.592 & 0.467 & 0.530 \\
SVR\_snp-06-004-052.mps & 0.088 & 0.246 & 0.488 \\
SVR\_snp-10-004-052.mps & 0.000 & 0.000 & 0.499 \\
SVR\_snp-10-052-052.mps & 0.000 & 0.487 & 0.549 \\
KNN\_snp-02-004-104.mps & 0.091 & 0.012 & 0.148 \\
KNN\_snp-04-052-052.mps & 0.310 & 0.467 & 0.530 \\
KNN\_snp-06-004-052.mps & 0.182 & 0.177 & 0.493 \\
KNN\_snp-10-004-052.mps & 0.004 & 0.000 & 0.492 \\
KNN\_snp-10-052-052.mps & 0.210 & 0.487 & 0.549 \\
RF\_snp-02-004-104.mps & 0.108 & 0.148 & 0.012 \\
RF\_snp-04-052-052.mps & 0.139 & 0.530 & 0.530 \\
RF\_snp-06-004-052.mps & 0.226 & 0.494 & 0.491 \\
RF\_snp-10-004-052.mps & 0.133 & 0.499 & 0.273 \\
RF\_snp-10-052-052.mps & 0.032 & 0.549 & 0.549 \\
MLP\_snp-02-004-104.mps & 0.138 & 0.012 & 0.148 \\
MLP\_snp-04-052-052.mps & 0.592 & 0.456 & 0.530 \\
MLP\_snp-06-004-052.mps & 0.088 & 0.191 & 0.490 \\
MLP\_snp-10-004-052.mps & 0.000 & 0.000 & 0.499 \\
MLP\_snp-10-052-052.mps & 0.000 & 0.549 & 0.549 \\
Stacking\_snp-02-004-104.mps & 0.108 & 0.148 & 0.389 \\
Stacking\_snp-04-052-052.mps & 0.179 & 0.530 & 0.799 \\
Stacking\_snp-06-004-052.mps & 0.191 & 0.415 & 0.491 \\
Stacking\_snp-10-004-052.mps & 0.133 & 0.499 & 0.499 \\
Stacking\_snp-10-052-052.mps & 0.032 & 0.549 & 0.565 \\
Voting\_snp-02-004-104.mps & 0.108 & 0.026 & 0.386 \\
Voting\_snp-04-052-052.mps & 0.592 & 0.530 & 0.821 \\
Voting\_snp-06-004-052.mps & 0.088 & 0.078 & 0.493 \\
Voting\_snp-10-004-052.mps & 0.036 & 0.000 & 0.499 \\
Voting\_snp-10-052-052.mps & 0.000 & 0.549 & 0.811 \\ \midrule
Average & \textbf{0.168} & 0.290 & 0.487 \\ \midrule
z-score &  & -4.743 & -12.965 \\
p-value &  & 1.05E-06 & 9.70E-39 \\ \bottomrule
\end{tabular}
\end{table}

\begin{landscape}
\begin{table}[]
\centering
\caption{Instance Features for PCA Analysis. All features, except for Density, Shape and Equality features were normalised using min-max normalisation on an instance by instance basis. Density and Equality features were not normalised, whilst the Shape features were normalised using min-max normalisation by considering all instances involved in the analysis.}
\begin{tabular}{ll}
\toprule
\multicolumn{1}{l}{\textbf{Features}} & \multicolumn{1}{l}{\textbf{Description}} \\ \hline
\toprule
Constr\_Sum\_Abs\_Obj\_mean & Mean sum of absolute values of objective coefficients associated with each constraint  \\
\hline 
Constr\_Sum\_Abs\_Obj\_stddev & Stddev sum of absolute values of objective associated with each constraint\\
\hline
Constr\_Sum\_Obj\_mean & Mean sum of values of objective coefficients associated with each constraint \\
\hline
Constr\_Sum\_Obj\_stddev & Standard Deviation sum of values of objective coefficients associated with each constraint \\
\hline
Obj\_terms\_mean & Mean of objective coefficients in the Instance \\
\hline
Obj\_terms\_stddev & Stddev of objective coefficients in the Instance \\
\hline
Non\_Zeroes\_mean & Mean No. Non Zeroes in constraints  \\
\hline
Non\_Zeroes\_stddev & Stddev of No. Non Zeroes in constraints \\ \hline
RHS\_Vals\_mean & Mean constraint RHS  \\ \hline
RHS\_Vals\_stddev & Stddev of constraint RHS \\ \hline
Density & No. Non Zeroes / (No. Variables * No. Constraints) \\ \hline
Shape & No. Variables / No. Constraints \\ \hline
Equality\_Prop & No. Equality Constraints / Total No. Constraints \\ \bottomrule
\end{tabular}
\label{tab: Instance Features}
\end{table}
\end{landscape}

\begin{landscape}
\begin{table*}[!htb]
\centering
\caption{Prediction RMSE Results: Shown for each Machine Learning based ranking method are the test RMSE scores found.}
\label{tab:RMSE Scores}
\begin{tabular}{@{}lllllllll@{}}
\toprule
                       & \multicolumn{8}{l}{\textbf{Ranking Method}}                                                                                      \\ \midrule
\textbf{Instance Name} & \textbf{Ridge} & \textbf{Lasso} & \textbf{SVR} & \textbf{KNN} & \textbf{RF} & \textbf{MLP} & \textbf{Stacking} & \textbf{Voting} \\ \midrule
cost266-UUE.mps        & 0.059          & 0.073          & 0.089        & 0.062        & 0.059       & 0.069        & 0.083             & 0.050           \\
dfn-bwin-DBE.mps       & 0.247          & 0.146          & 0.128        & 0.189        & 0.118       & 0.123        & 0.107             & 0.140           \\
germany50-UUM.mps      & 0.144          & 0.073          & 0.113        & 0.083        & 0.088       & 0.099        & 0.084             & 0.100           \\
ta1-UUM.mps            & 0.244          & 0.082          & 0.153        & 0.126        & 0.097       & 0.132        & 0.090             & 0.137           \\
ta2-UUE.mps            & 0.080          & 0.046          & 0.065        & 0.047        & 0.044       & 0.030        & 0.039             & 0.036           \\ \midrule
g200x740.mps           & 0.113          & 0.119          & 0.084        & 0.110        & 0.123       & 0.101        & 0.181             & 0.070           \\
h50x2450.mps           & 0.197          & 0.059          & 0.054        & 0.227        & 0.064       & 0.041        & 0.096             & 0.076           \\
h80x6320d.mps          & 0.245          & 0.356          & 0.278        & 0.206        & 0.302       & 0.243        & 0.335             & 0.254           \\
k16x240b.mps           & 0.157          & 0.143          & 0.134        & 0.116        & 0.099       & 0.077        & 0.078             & 0.090           \\ \midrule
snp-02-004-104.mps     & 0.239          & 0.203          & 0.241        & 0.217        & 0.207       & 0.222        & 0.216             & 0.219           \\
snp-04-052-052.mps     & 0.238          & 0.181          & 0.160        & 0.144        & 0.170       & 0.154        & 0.147             & 0.155           \\
snp-06-004-052.mps     & 0.137          & 0.148          & 0.161        & 0.133        & 0.115       & 0.156        & 0.104             & 0.119           \\
snp-10-004-052.mps     & 0.135          & 0.161          & 0.138        & 0.129        & 0.100       & 0.135        & 0.113             & 0.116           \\
snp-10-052-052.mps     & 0.187          & 0.196          & 0.156        & 0.174        & 0.193       & 0.161        & 0.173             & 0.159           \\ \bottomrule
\end{tabular}
\end{table*}
\end{landscape}

\makeatletter\@input{xx.tex}\makeatother